# On the Fine Structure of Quadratical Quasigroups

## R. A. R. Monzo[1]

**Abstract.** We prove that quadratical quasigroups form a variety $Q$ of right and left simple groupoids. New examples of quadratical quasigroups of order 25 and 29 are given. General observations of the fine structure of quadratical quasigroups are described and inter-relationships between their properties are explored. The spectrum of $Q$ is proved to be contained in the set of integers equal to 1 plus a multiple of 4.

## 1. Introduction

This paper builds on the work of Polonijo [2], Volenec [3] and Dudek [1] on quadratical quasigroups. Polonijo [2] and Volenec [3] proved that a quadratical groupoid is a quasigroup. Volenec [3,4] gave a motivation for studying quadratical quasigroups, in terms of a geometrical representation of the complex numbers $C$ as points of the Euclidean plane. He defined a product $*$ on $C$ that defines a quadratical quasigroup and in which the product of distinct elements $x$ and $y$ is the third vertex of a positively oriented, isosceles right triangle, at which the right angle occurs. He proved a number of properties of quadratical quasigroups, which are listed in Theorem 2.1 below. These properties tell us a great deal and, indeed, we apply them to prove that quadratical quasigroups form a variety $Q$ (Corollaries 2.4 and 2.5 and Theorem 2.11). Inter-relationships amongst the properties of quadratical quasigroups are explored in Section 2.

We also amplify our understanding of the fine structure of quadratical quasigroups. In so doing, we give further meaning to the "quad" in the word quadratical, in terms of "4-cycles". We apply this to prove that the order of a finite quadratical quasigroup is $1+4n$ for some $n \in \{0,1,2,...\}$ (Propositions 3.1-3.4), fine tuning Dudek's result that the order of a finite quadratical quasigroup is odd [1, Corollary 1].

Dudek [1] determined all quadratical groupoids of order less than 24. We prove in the remainder of the paper that the quadratical quasigroups of orders 5, 9, 13, and 17 discovered by Dudek are quadratical quasigroups of form $Qn$ (defined below), for some $n \in \{1,2,3,4\}$, or their dual groupoids. We give several examples of quadratical quasigroups of order 25 and introduce the concept of a translatable groupoid, applying it to find a quadratical quasigroup of order 29. The paper ends with some open questions, a conjecture about the spectrum of $Q$ and ideas about the future direction of research.

## 2. Inter-relationships between properties of quadratical groupoids.

**Definition 2.1.** A groupoid $(G, \cdot)$ has *property A* if it satisfies the identity $xy \cdot x = zx \cdot yz$. It is called *right solvable* [*left solvable*] if for any $\{a,b\} \subseteq G$ there exists a *unique* $x \in G$ such that $ax = b$ [ $xa = b$ ]. It is *left* [*right*] *cancellative* if $xy = xz$ implies $y = z$ [ $yx = zx$ implies $y = z$ ]. It is a *quasigroup* if it is left and right solvable.

Note that a right solvable groupoid is left cancellative and a left solvable groupoid is right cancellative.

Volenec [3] defined a quadratical groupoid as a right solvable groupoid satisfying property A. He proved that *a quadratical groupoid is left solvable* and satisfies the following identities:

---

**Theorem 2.1 [3, Theorems 1,2,3,4 and 5]:** *A quadtratical groupoid is left solvable and satisfies the following identities:*

(1) $x = x^2$ (*idempotency*);
(2) $x \cdot yx = xy \cdot x$ (*elasticity*);
(3) $x \cdot yx = xy \cdot x = yx \cdot y$ (*strong elasticity*);
(4) $yx \cdot xy = x$ (*bookend*);

(5) $x \cdot yz = xy \cdot xz$ (*left distributivity*);
(6) $xy \cdot z = xz \cdot yz$ (*right distributivity*);
(7) $xy \cdot zw = xz \cdot yw$ (*mediality*);
(8) $x(y \cdot yx) = (xy \cdot x)y$;
(9) $(xy \cdot y)x = y(x \cdot yx)$ and
(10) $xy = zw$ if and only if $yz = wx$ (*alterability*).

**Corollary 2.2. [2 and 3, Theorem 5]** *A quadratical groupoid is a quasigroup.*

Note that throughout the remainder of this paper we will use the fact that quadratical groupoids are quasigroups and satisfy properties (1) through (10), often without mention. Note also that property (3) allows us to write the term *xyx* without ambivalence in any quadratical quasigroup.

**Definition 2.2.** We define $Q$ to be the collection of ***quadratical quasigroups***.

**Theorem 2.3.** Q *is a quadratical quasigroup if and only if* Q *is a bookend, medial, strongly elastic groupoid satisfying property A.*

Proof: ($\Rightarrow$) This follows from Theorem 2.1 and the definition of a quadratical quasigroup.

($\Leftarrow$) First we prove that Q is left cancellative. Suppose that $ax = ay$. Then,

$x = ax \cdot xa$    (bookend)
$\phantom{x} = ay \cdot xa$    ($ax = ay$)
$\phantom{x} = yx \cdot y$    (property A)
$\phantom{x} = xy \cdot x$    (strong elasticity)
$\phantom{x} = ax \cdot ya$    (property A)
$\phantom{x} = ay \cdot ya$    $(ax = ay)$
$\phantom{x} = y$       (bookend).
So $x = y$ and Q is left cancellative.

Property A implies $x^2 \cdot x = x^2 x^2$ and so left cancellativity implies $x = x^2$. Hence, Q is idempotent. Since Q is medial and idempotent it is therefore (left and right) distributive.

Using left and right distributivity, mediality, idempotency and strong elasticity we have
$a(b \cdot ba) = (ab)(a \cdot ba) = (ab)(ba \cdot b) = (a \cdot ba)b = (ab \cdot a)b$. Hence,

$$(1): a(b \cdot ba) = (ab \cdot a)b.$$

We now prove that $ax = b$ has a unique solution $x = (b \cdot ba) \cdot (b \cdot ba)(ba \cdot a)$.

$ax = a(b \cdot ba) \cdot [a(b \cdot ba) \cdot a(ba \cdot a)]$   (distributivity)

$= (aba \cdot b)[(aba \cdot b)(aba \cdot a)]$   (from (1), distributivity and strong elasticity)

$= (aba)(b \cdot ba)$   (from distributivity)

$= ab \cdot ba$   (from strong elasticity, mediality and idempotency)

$= b$   (bookend)

The solution $x$ is unique because Q is left cancellative. We have proved that Q is right solvable and so by definition, Q is a quadratical quasigroup. ∎

**Corollary 2.4.** Q *is a quadratical quasigroup if and only if* Q *is a bookend, medial, strongly elastic groupoid that satisfies property A.*

**Corollary 2.5.** Q *is a quadratical quasigroup if and only if* Q *is a medial, idempotent groupoid that satisfies property A.*

Proof: $(\Rightarrow)$ This follows from the definition of a quadratical quasigroup and from Theorem 2.1.

$(\Leftarrow)$ Let Q be a medial, idempotent groupoid that satisfies property A. By Corollary 2.4 we need only show that Q is strongly elastic and bookend.

Since Q is idempotent and satisfies property A, for all $\{x, y\} \subseteq Q$, $x = x^2 = x^2 x = yx \cdot xy$ and Q is bookend.

Also, idempotent and medial implies distributive and so $xy \cdot x = xy \cdot x^2 = x^2 \cdot yx = x \cdot yx = yx \cdot yy = yx \cdot y$ and Q is strongly elastic. ∎

**Theorem 2.6.** *A distributive, bookend groupoid* G *is idempotent*

Proof: Since G is bookend, $xx \cdot xx = x$ for any $x \in G$. Since G is distributive,
$x^2 x = x^2(x^2 x^2) = (x^2 x^2)(x^2 x^2) = x^2 = (x^2 x)x = (x^2 x)x^2 = x^2 x^2 = x$. So G is idempotent. ∎

**Example 2.1.** A groupoid G of order greater than or equal to 2 and satisfying the identity $xy = zw$ ($x, y, z, w \in G$) is distributive but not idempotent.

**Example 2.2.** The following groupoid G is bookend but is *not* idempotent, elastic, medial, distributive, alterable, left solvable or right solvable.

| G | x | y | z | w |
|---|---|---|---|---|
| x | y | z | w | y |
| y | w | x | w | x |
| z | y | x | w | y |
| w | z | z | x | z |

**Theorem 2.7.** *A distributive, medial, bookend groupoid* G *is cancellative.*

Proof: Suppose that $ax = ay$ for any $\{a, x, y\} \subseteq G$. By Theorem 2.6, $ax = (ax)^2$. Then,
$ax = ax \cdot ax = ax \cdot ay = ay \cdot ax = a(yx) = a(xy)$. Then, $yx = (a \cdot yx)(yx \cdot a) = (ax)(yx \cdot a) = (ax \cdot yx)(ax \cdot a) =$

$= (ay \cdot x)(ay \cdot a) = (ay)^2(xa) = ay \cdot xa = ax \cdot xa = x$. Similarly, $xy = y$. So $y = xy \cdot yx = yx = x$.
Similarly, $xa = ya$ for any $\{a, x, y\} \subseteq G$ implies $x = y$, so G is cancellative. ∎

**Theorem 2.8.** *A distributive, bookend groupoid G is strongly elastic.*

Proof: Theorems 2.6 and 2.7 imply that G is idempotent and cancellative. So $xy \cdot x = xy \cdot xx = x \cdot yx$. Hence, $(xy \cdot x)y = (x \cdot yx)y = (xy)(yx \cdot y) = (xy \cdot yx)(xy \cdot y) = y(xy \cdot y) = (y \cdot xy)y^2 = (yx \cdot y)y$ and, by cancellation, $xy \cdot x = yx \cdot y$. Therefore, G is strongly elastic.

**Corollary 2.9.** *Medial, idempotent, bookend groupoids are cancellative and strongly elastic.*

Proof: Medial idempotent groupoids are distributive. The corollary follows from Theorems 2.7 and 2.8. ∎

**Theorem 2.10.** *A left (or right) cancellative, medial, idempotent, strongly elastic groupoid is bookend.*

Proof: Mediality and idempotency imply distributivity. Then, $(ca \cdot ac)^2 = ca \cdot ac = (ca \cdot a)(ca \cdot c)$. Then strong elasticity implies $(ca \cdot ac)^2 = ca \cdot ac = (ca \cdot a)(c \cdot ac) = (ca \cdot a)(ca \cdot c) = (ca \cdot a)(ac \cdot a) = (ca \cdot ac)a$ and so left cancellativity implies $ca \cdot ac = a$. Also, using right cancellativity,
$(ca \cdot ac)^2 = ca \cdot ac = ca \cdot (ac)^2 = (c \cdot ac)(a \cdot ac) = (a \cdot ca)(a \cdot ac) = a(ca \cdot ac)$ implies $a = ca \cdot ac$. ∎

**Theorem 2.11.** Q *is a quadratical quasigroup if and only if* Q *is idempotent, bookend and medial.*

Proof: $(\Rightarrow)$ This follows from Theorem 2.1.

$(\Leftarrow)$ By Corollary 2.5, we need only show that $xy \cdot x = zx \cdot yz$. But, since medial and idempotent implies distributive, $zx \cdot yz = (zx \cdot y)(zx \cdot z)$ and, by Corollary 2.9, $zx \cdot yz = (zx \cdot y)(xz \cdot x) = (zx \cdot xz)(yx) = x \cdot yx = xy \cdot x$. ∎

**Definition 2.3:** *The dual of a groupoid* $(G, \cdot)$ is the groupoid $G^* = (G, *)$ where $x * y = x \cdot y$ $(x, y \in G)$. A groupoid will be called *self-dual* if it is isomorphic to its dual groupoid.

The dual of an idempotent [bookend; medial] groupoid is idempotent [bookend; medial]. Hence,

**Corollary 2.12.** *The dual of a quadratical quasigroup is quadratical.*

**Corollary 2.13.** *Any subgroupoid of a quadratical quasigroup is quadratical.*

Note that a semigroup is idempotent and bookend if and only if it satisfies the identity $x = xyx$; that is, if and only if it is a rectangular band. A semigroup is cancellative, with property A, if and only if it is trivial.

**Theorem 2.14.** *Idempotent, bookend, alterable groupoids are elastic.*

Proof: $x = x^2 = yx \cdot xy$ implies $x \cdot yx = xy \cdot x$. ∎

**Theorem 2.15.** *Elastic bookend groupoids are idempotent.*

Proof: $x^2 = xx^2 \cdot x^2 x = x^2 x \cdot xx^2 = x$. ∎

**Theorem 2.16.** Q *is a quadratical quasigroup if and only if Q is elastic, bookend and medial.*

Proof: $(\Rightarrow)$ This follows from Theorem 2.1.
$(\Leftarrow)$ Assume that Q is elastic, bookend and medial. By Theorem 2.15, Q is idempotent. By Theorem 2.11 then, Q is quadratical. ∎

**Theorem 2.17.** *Elastic, idempotent, alterable groupoids are bookend.*

Proof: $x \cdot yx = xy \cdot x$ implies $yx \cdot xy = x^2 = x$. ∎

**Corollary 2.18.** *If* G *is an idempotent, alterable groupoid then* G *is elastic if and only if it is bookend.*

Proof: This follows from Theorems 2.14 and 2.17. ∎

**Theorem 2.19.** *Medial bookend groupoids are alterable.*

Proof: Suppose that $xy = zw$. Then $zw \cdot yx = xy \cdot yx = y$ and $wz \cdot xy = wz \cdot zw = z$. Therefore, $yz = (zw \cdot yx)(wz \cdot xy) = (zw \cdot wz)(yx \cdot xy) = wx$. ∎

**Theorem 2.20.** *A groupoid* Q *is a quadratical quasigroup if and only if it is elastic, medial, idempotent and alterable.*

Proof: $(\Rightarrow)$ This follows from Theorem 2.1.
$(\Leftarrow)$ Suppose that Q is elastic, medial, idempotent and alterable. By Theorem 2.17, Q is bookend. By Theorem 2.16, Q is quadratical. ∎

**Theorem 2.21.** *A left (or right) distributive, alterable, strongly elastic groupoid satisfies property A.*

Proof: Since in a left distributive groupoid $y \cdot zx = yz \cdot yx$, $zx \cdot yz = yx \cdot y = xy \cdot x$. Similarly, a right distributive, alterable, strongly elastic groupoid satisfies property A. ∎

**Theorem 2.22.** *In a quadratical quasigroup* $x \cdot yz = xy \cdot z$ *if and only if* $x = z$.

Proof: $(\Rightarrow)$ Using Theorem 2.1, $x \cdot yz = xy \cdot z$ implies $xy \cdot xz = xz \cdot yz$ implies $yz \cdot xy = (xz)^2 = xz = zx \cdot z$ implies $xz = zx \cdot z$ implies $x = zx$ implies $x = z$.
$(\Leftarrow)$ By Theorem 2.1, a quadratical quasigroup is elastic and so $x \cdot yx = xy \cdot x$. ∎

**Definition 2.4.** A groupoid is ***nowhere commutative*** if $xy = yx$ implies $x = y$.

**Theorem 2.23.** *Quadratical quasigroups are nowhere commutative.*

Proof: Since by Theorem 2.1, quadratical quasigroups are alterable and idempotent, $xy = yx$ implies $y^2 = x^2$ implies $y = x$. ∎

**Theorem 2.24.** Q *is a quadratical quasigroup if and only if Q is distributive, bookend and alterable.*

Proof: $(\Rightarrow)$ This follows from Theorem 2.1.

$(\Leftarrow)$ By Theorem 2.6, Q is idempotent. Therefore, by Theorem 2.11 we need only show that Q is medial. Now distributive bookend groupoids are strongly elastic by Theorem 2.8. By Theorem 2.21, left distributive, alterable, strongly elastic groupoids satisfy property A. So Q satisfies property A. Hence, $wx \cdot w = zw \cdot xz = yw \cdot xy$ and, using alterability, $xz \cdot yw = xy \cdot zw$. So Q is medial. ∎

**Definition 2.5.** A subset *I* of a groupoid G is a ***right [left] ideal of* G** if $ig \in I$ $[gi \in I]$ for all $i \in I$ and all $g \in G$. The subset *I* is called an ideal of G if it is a right ideal and a left ideal of G. *A groupoid **G** is called simple [right simple; left simple]* if for every ideal [right ideal; left ideal] *I* of G, $I = G$.

**Theorem 1.25.** *Bookend groupoids are right simple and left simple groupoids.*

Proof: Suppose that *I* is a right or left ideal of a bookend groupoid G. Let $i \in I$ and $g \in G$. Then, $g = ig \cdot gi \in I$ and so $I = G$. ∎

**Corollary 1.26.** *Quadratical quasigroups are right and left simple.*

### 3. The fine structure of quadratical quasigroups.

Let Q be a quadratical quasigroup with $\{a,b\} \subseteq Q$ and $a \neq b$. Suppose that $C = \{x_1, x_2, ..., x_n\} \subseteq Q$ consists of $n$ distinct elements, such that $aba = x_1 \cdot x_2 = x_2 \cdot x_3 = x_3 \cdot x_4 = ... = x_{n-1} \cdot x_n = x_n \cdot x_1$. Then $C$ will be called an (ordered) *n*-cycle based on *aba*. Note that $x_1 \neq aba$, or else $x_1 = x_2 = x_3 = x_4 = ... = x_n = aba$. Note also that if $C = \{x_1, x_2, ..., x_n\} \subseteq Q$ is an *n*-cycle based on *aba* then so is $Ci = \{x_i, x_{(i+1, \bmod n)}, x_{(i+2, \bmod n)}, ..., x_{(i+n-1, \bmod n)}\}$.

**Proposition 3.1.** *If n-cycles exist in a quadratical quasigroup then* $n = 4$.

Proof. Since $aba = x_n \cdot x_1 = x_1 \cdot x_2 = x_2 \cdot x_3$, by alterability $x_1 = x_2 \cdot x_n$ and $x_2 = x_3 \cdot x_1$. Now $x_3(x_2 \cdot x_4) = (x_3 \cdot x_2)(x_3 \cdot x_4) = (x_3 \cdot x_2) \cdot aba = (x_3 \cdot aba)(x_2 \cdot aba)$. But by alterability, $aba \cdot x_2 = x_3 \cdot aba$ and so $x_3(x_2 \cdot x_4) = (x_3 \cdot aba)(x_2 \cdot aba) = (aba \cdot x_2)(x_2 \cdot aba) = x_2 = x_3 \cdot x_1$. Hence, by cancellation, $x_1 = x_2 \cdot x_4 = x_2 \cdot x_n$ and so $x_4 = x_n$. ∎

**Proposition 3.2.** *Let* Q *be a quadratical quasigroup with* $\{a,b\} \subseteq Q$ *and* $a \neq b$. *Then every element* $x_1 \neq aba$ *of* Q *is a member of a 4-cycle based on aba.*

Proof. Let $\{a,b\} \subseteq Q$ and $a \neq b$. Suppose that $x_1 \neq aba$ for some $x_1 \in Q$. Using right solvability, we can solve the equations $aba = x_1 \cdot x$, $aba = xy$, $aba = yz$ and $aba = zw$. If we define $x_2 = x$, $x_3 = y$, $x_4 = z$ and $x_5 = w$ then $aba = x_1 \cdot x_2 = x_2 \cdot x_3 = x_3 \cdot x_4 = x_4 \cdot x_5$. Using alterability, $x_4 = x_5 x_3$ and $x_5 x_1 = x_2 x_4 = x_2(x_5 x_3) = (x_2 x_5)(x_2 x_3) = (x_2 x_5) aba$. Therefore, $aba \cdot x_5 = x_1(x_2 x_5) = (x_1 x_2)(x_1 x_5) = aba \cdot (x_1 x_5)$. Hence $x_5 = x_1 x_5$ and $x_1 = x_5$. So we have proved that $\{x_1, x_2, x_3, x_4\}$ is a 4-cycle based on $aba$. ∎

**Proposition 3.3.** *Let C and D be two 4-cycles based on aba ($a \neq b$) in a quadratical quasigroup. Then either $C = D$ or $C \cap D = \varnothing$.*

Proof. Suppose that $C = \{x_1, x_2, x_3, x_4\}$ and $D = \{y_1, y_2, y_3, y_4\}$. If $x_1 = y_1$ then $aba = x_1 \cdot x_2 = y_1 \cdot y_2 = x_1 \cdot y_2$ and so $x_2 = y_2$. Then, $aba = x_2 \cdot x_3 = y_2 \cdot x_3 = y_2 \cdot y_3$ and so $x_3 = y_3$. Finally, $aba = x_3 \cdot x_4 = y_3 \cdot x_4 = y_3 \cdot y_4$ and $x_4 = y_4$. Hence, $C = D$. Similarly, if $x_1 = y_2$ then we can prove that $x_2 = y_3$, $x_3 = y_4$ and $x_4 = y_1$ and $C = D$. Similarly, if $x_1 \in \{y_3, y_4\}$ it is straightforward to prove that $C = D$. The proofs that $C = D$ if $x_2 \in D$ or $x_3 \in D$ or $x_4 \in D$ are similar. ∎

**Proposition 3.4.** *A finite quadratical quasigroup has order $1+4n$ for some $n \in \{0,1,2,...\}$.*

Proof. A finite quadratical quasigroup consists of the element $aba$ and the union of its disjoint 4-cycles based on $aba$. By definition, no cycle contains the element $aba$. The proposition is therefore valid. ∎

**Definition 3.1.** Let Q be a quadratical quasigroup with $\{a,b\} \subseteq Q$ and $a \neq b$. Then $\{a, b, ab, ba, aba\}$ contains five distinct elements. We use the notation **[1,1] = $a$**, **[1,2] = $ab$**, **[1,3] = $ba$** and **[1,4] = $b$**. We omit the commas and square brackets in the notation, when this causes no confusion, and write $11 = a$, $12 = ab$, $13 = ba$ and $14 = b$. For $n \geq 2$, by induction we define **n1 = (n-1)1 · (n-1)2, n2 = (n-1)2 · (n-1)4, n3 = (n-1)3 · (n-1)1 and n4 = (n-1)4 · (n-1)3** and **H$n$ = {$n$1, $n$2, $n$3, $n$4}**. On the occasions when we need to highlight that the element $fk$ ($f \in \{1,2,...n\}$ and $k \in \{1,2,3,4\}$) is in the dual quadratical quasigroup $Q = (Q,*)$ we will denote it by $fk^*$. Similarly, **H$n$* = {$n$1*, $n$2*, $n$3*, $n$4*}**. Note that the values of both $fk$ and $fk^*$ depend on the choice of the elements $a$ and $b$, as in Examples 3.1 and 3.2.

**Example 3.1.** $H2 = \{a \cdot ab,\ ab \cdot b,\ ba \cdot a,\ b \cdot ba\}$,
$H3 = \{(a \cdot ab)(ab \cdot b),\ (ab \cdot b)(b \cdot ba),\ (ba \cdot a)(a \cdot ab),\ (b \cdot ba)(ba \cdot a)\}$ and
$H4 = \{[(31)(32) \cdot (32)(34),\ (32)(34) \cdot (34)(33),\ (33)(31) \cdot (31)(32),\ (34)(33) \cdot (33)(31)]\}$, where
$31 = (a \cdot ab)(ab.b)$, $32 = (ab \cdot b)(b \cdot ba)$, $33 = (ba \cdot a)(a \cdot ab)$ and $34 = (b \cdot ba)(ba \cdot a)$.

**Example 3.2.** $11^* = a$, $12^* = a*b$, $13^* = b*a$, $14^* = b$ and, for $n \geq 2$, by induction we define
$n1^* = (n-1)1^* \ * \ (n-1)2^*$, $n2^* = (n-1)2^* \ * \ (n-1)4^*$, $n3^* = (n-1)3^* \ * \ (n-1)1^*$ and $n4^* = (n-1)4^* \ * \ (n-1)3^*$.

**Example 3.3.** $H2^* = \{a*(a*b),\ (a*b)*b,\ (b*a)*a,\ b*(b*a)\} = \{ba \cdot a,\ b \cdot ba,\ a \cdot ab,\ ab \cdot b\}$ and
$52^* = 42^* 44^* = (32^* 34^*)(34^* 33^*) =$
$= \{ [(ab^* b)^* (b^* ba)] \ * \ [(b^* ba)^* (ba^* a)] \} \ * \ \{ [(b^* ba)^* (ba^* a)] \ * \ [(ba^* a)^* (a^* ab)] \}$, where
$a^* ab = a*(a*b)$, $ab^* b = (a*b)*b$, $ba^* a = (b*a)*a$ and $b^* ba = b*(b*a)$.

Note that the expression *ab*, when working in the dual groupoid Q = (Q,*) , equals $a*b$, which equals $b \cdot a$ in the original groupoid itself. This notation will cause no problems, as we will either calculate values only using the dot product or the star product, or when we are calculating using both products, as in Theorem 5.1, the distinction will be obvious.

The proofs of the following propositions are straightforward, using induction on *n* and the properties of quadratical quasigroups, and are omitted.

**Proposition 3.5.** *For any positive integer n, $n1 \cdot n4 = n2$, $n2 \cdot n3 = n4$, $n3 \cdot n2 = n1$ and $n4 \cdot n1 = n3$.*

**Proposition 3.6.** *For $n > 1$, $aba \cdot nk = (n-1)k$ for any $k \in \{1,2,3,4\}$.*

**Proposition 3.7.** *For $n > 1$, $n1 \cdot aba = (n-1)2$, $n2 \cdot aba = (n-1)4$, $n3 \cdot aba = (n-1)1$ and $n4 \cdot aba = (n-1)3$..*

**Proposition 3.8.** *For any positive integer n, Hn contains 4 distinct elements.*

**Proposition 3.10.** *For any positive integer n, $Hn \cap \{aba\} = \emptyset$.*

**Proposition 3.11.** *For any positive integer n, $n1 \cdot n3 = n2 \cdot n1 = n3 \cdot n4 = n4 \cdot n2 = aba$.*

**Proposition 3.12.** *$Hn = \{n1, n3, n4, n2\}$ is a 4-cycle based on aba.*

**Definition.** We say that *a quadratical quasigroup* **Q** *is of the form Qn*, for some positive integer *n*, if $Q = \{aba\} \cup \bigcup_{i=1}^{n} Hi$ for some $\{a,b\} \subseteq Q$.

Dudek [1] proved that there is no quadratical quasigroup of order 21 and so there is no quadratical quasigroup of form $Q5$. We show in the next section that the quadratical quasigroups of orders 5, 9, 13 and 17 found by Dudek [1], and their dual groupoids, are of the form $Q1$, $Q2$, $Q3$ and $Q4$ respectively.

## 4. Products in quadratical quasigroups.

We are now in a position to examine more closely the Cayley tables of quadratical quasigroups. This will aid in the construction of the tables for quadratical quasigroups of orders 5, 9, 13 and 17. Dudek [1] determined all quadratical quasigroups of order $\leq 24$. He gave two examples of quadratical quasigroups of orders 5, 13 and 17 and six examples of quadratical quasigroups of order 9. A close examination of the fine structure of quadratical quasigroups will aid us in proving that the quadratical quasigroups of Dudek are all of the form $Qn$, for some positive integer *n*. Each pair of quadratical quasigroups of orders 5, 13 or 17 will be proved to be dual groupoids. The 6 quadratical quasigroups of order 9 will be proved to be of form $Q2$ and self-dual. That is, up to isomorphism, there is only one quadratical quasigroup of order 9.

A method of constructing quadratical quasigroups of the form $Qn$ is as follows. Proposition 3.6 implies that $aba \cdot Hn = H(n-1)$ for all $n \neq 1$. Since quadratical quasigroups are cancellative, we can assume that $aba \cdot H1 = Hn$. If we choose the value of $aba \cdot 11$ in $Hn = \{n1, n2, n3, n4\}$ then, using the properties of quadratical quasigroups, we can attempt to fill in the remaining unknown products in the Cayley table. If this can be done without contradiction, then (using Theorem 2.11) we can check that the groupoid thus obtained is quadratical, by checking that it is bookend and medial. Completing the Cayley table is this way is **_not_** always possible, as shown in the following example.

**Example 4.1.** Suppose Q2 is a quadratical quasigroup of form $Q2$. Then
$aba \cdot 11 = aba \cdot a \in H2 = \{21, 22, 23, 24\} = \{a \cdot ab, ab \cdot b, ba \cdot a, b \cdot ba\}$. Now $aba \cdot a = a(ba \cdot a)$ and so $aba \cdot a \notin \{a \cdot ab, ba \cdot a\}$, since cancellativity, idempotency and alterability would imply that $a = b$ (if $aba \cdot a = ba \cdot a$) and $b = a \cdot ab$ (if $aba \cdot a = a \cdot ab$), the latter contradicting Proposition 3.3. Hence, $aba \cdot a$ must be in the set $\{ab \cdot b, b \cdot ba\}$. However, If $aba \cdot a = b \cdot ba$ then by alterability,
$ab = ba \cdot aba = (b \cdot ab)a = aba \cdot a = b \cdot ba$, a contradiction since $H1 \cap H2 = \varnothing$ in Q2, if $Q2$ is quadratical.

Example 4.1 shows that $aba \cdot a = ab \cdot b$. Using the properties of quadratical quasigroups, the Cayley table of the groupoid we will call Q2 can only be completed in one way, as shown below here.

**Table 1.**

| Q2 | 11 = a | 12 = ab | 13 = ba | 14 = b | aba | 21 = a·ab | 22 = ab·b | 23 = ba·a | 24 = b·ba |
|---|---|---|---|---|---|---|---|---|---|
| 11 = a | a | a·ab | aba | ab | b·ba | ba | b | ab·b | ba·a |
| 12 = ab | aba | ab | b | ab·b | ba·a | b·ba | a | a·ab | ba |
| 13 = ba | ba·a | a | ba | aba | ab·b | ab | b·ba | b | a·ab |
| 14 = b | ba | aba | b·ba | b | a·ab | ab·b | ba·a | a | ab |
| aba | ab·b | b·ba | a·ab | ba·a | aba | a | ab | ba | b |
| 21 = a·ab | b·ba | b | ba·a | a | ab | a·ab | ba | aba | ab·b |
| 22 = ab·b | a·ab | ba·a | ab | ba | b | aba | ab·b | b·ba | a |
| 23 = ba·a | ab | ba | ab·b | b·ba | a | b | a·ab | ba·a | aba |
| 24 = b·ba | b | ab·b | a | a·ab | ba | ba·a | aba | ab | b·ba |

We then need to calculate all the possible products $xy \cdot yx$ and $xy \cdot zw$ in Table 1, to prove that they are equal to $y$ and $xz \cdot yw$ respectively. Then, by Theorem 2.11, Q2 would be quadratical. This proves to be the case and we omit the detailed calculations. However, to give a flavour of the calculations we find all products $aba \cdot x$ and $x \cdot aba$ when $x \in H1$ and $aba \cdot a = ab \cdot b$.

Since $(a \cdot aba)(aba \cdot a) = (a \cdot aba)(ab \cdot b)$, it follows that $a \cdot aba = b \cdot ba$, $aba \cdot b = ba \cdot a$, $aba \cdot ab = (aba \cdot a)(aba \cdot b) = (ab \cdot b)(ba \cdot a) = b \cdot ba$ and, similarly $aba \cdot ba = a \cdot ab$. Then $aba \cdot ab = b \cdot ba$ implies $ba \cdot aba = ab \cdot b$. Also, $aba = (ab \cdot aba)(aba \cdot ab) = (ab \cdot aba)(b \cdot ba)$ implies $ab \cdot aba = ba \cdot a$. Finally, $b \cdot aba = (ab \cdot aba)(ba \cdot aba) = (ba \cdot a)(ab \cdot b) = a \cdot ab$.

**Open question.** Assume that the Cayley table for a quadratical quasigroup Q of form $Qn$ can be completed, without contradiction, using the properties of quadratical quasigroups. We speculate that it should follow that Q is quadratical. We cannot prove this but we conjecture that a proof exists.

Now, if we calculate the Cayley table for $(Q2)^*$, the dual of Q2, we see that the table for the dual product $*$ (defined as $a * b = b \cdot a$) is exactly the same as Table 1, where the product is the dot product $\cdot$. (For example, $[(b*a)*a]*(b*a) = [a \cdot (a \cdot b)]*(a \cdot b) = (a \cdot b) \cdot [a \cdot (a \cdot b)] = [b \cdot (b \cdot a)] = (a*b)*b$ and, by Table 1, $[(b \cdot a) \cdot a] \cdot (b \cdot a) = [(a \cdot b) \cdot b]$). Hence, $Q2 \cong (Q2)^*$. Another way to put this is that the quadratical quasigroup $Q2$ must be self-dual. An isomorphism $\theta$ between $Q2$ and $(Q2)^*$ is: $\theta a = a$, $\theta b = b$, $\theta(ab) = a*b$, $\theta(ba) = b*a$, $\theta(a \cdot ab) = a*(a*b)$, $\theta(ab \cdot b) = (a*b)*b$, $\theta(ba \cdot a) = (b*a)*a$ and $\theta(b \cdot ba) = b*(b*a)$.

**Example 4.2.** It is straightforward to calculate the Cayley tables of the quadratical quasigroups, each of order 9, given by Dudek [1]. They are each based on the group $Z_3 \times Z_3$ of ordered pairs of integers, with product being addition (mod 3). The products are defined as follows:

$(x,y)*(z,u) = (y+z+2u, \ x+y+2z)$, $\qquad (x,y)*(z,u) = (2y+z+u, \ 2x+y+z)$,

$(x,y)*(z,u) = (x+y+2u, \ x+2z+u)$, $\qquad (x,y)*(z,u) = (x+2y+u, \ 2x+z+u)$,

$(x,y)*(z,u) = (2x+y+2z+2u, \ 2x+2y+z+2u)$ and $(x,y)*(z,u) = (2x+2y+2z+u, \ x+2y+2z+2u)$.

In each table, if we calculate $ab$ and $ba$ for the ordered pairs $a = (1,1)$ and $b = (1,2)$ we see that $Q = \{aba\} \cup H1 \cup H2$ and that $aba \cdot a = ab \cdot b$. Therefore, these six quadratical quasigroups are isomorphic to each other and to $Q2$. We already knew that there is only one quadratical quasigroup of order 9, but these calculations clarify (and reinforce a conviction) that the quadratical quasigroups of order 9 found by Dudek are isomorphic.

**Example 4.3.** We now calculate the Cayley table for a groupoid we call Q1 and its dual, when $aba \cdot a \in \{ab \cdot b\}$.

**Table 2.**

| Q1  | a   | ab  | ba  | b   | aba |     |     | (Q1)* | a   | a*b | b*a | b   | aba |
|-----|-----|-----|-----|-----|-----|-----|-----|-------|-----|-----|-----|-----|-----|
| a   | a   | ba  | aba | ab  | b   | ……  | ……  | a     | a   | b   | aba | a*b | b*a |
| ab  | aba | ab  | b   | a   | ba  | ……  | ……  | a*b   | aba | a*b | b   | b*a | a   |
| ba  | b   | a   | ba  | aba | ab  | ……  | ……  | b*a   | a*b | a   | b*a | aba | b   |
| b   | ba  | aba | ab  | b   | a   | ……  | ……  | b     | b*a | aba | a   | b   | a*b |
| aba | ab  | b   | a   | ba  | aba | ……  | ……  | aba   | b   | b*a | a*b | a   | aba |

Checking these tables shows that each is medial and bookend and that, indeed, these two quadratical quasigroups are dual.

**Open question.** Examining Tables 1 and 2 closely, we can show that any two distinct elements of Q1 $[(Q1)^*; Q2]$ generate Q1 $[(Q1)^*; Q2]$. This will later be seen to be the case also for $Q3, Q4$ and their duals. We conjecture that if $Q$ is a quadratical quasigroup of form $Qn$, for some positive integer $n$, then it is generated by any two distinct elements. Such a property does **not** hold in quadratical quasigroups in general, as we shall now prove.

**Example 4.4.** Since $Q$ is a variety of groupoids, the direct product of quadratical quasigroups is quadratical. Hence, Q1×Q1 is quadratical. If we choose a base element, $(a,b)$ say, then Q1×Q1 consists of six disjoint 4-cycles based on $(a,b)$; namely, $\{(a,a),(a,aba),(a,ab),(a,ba)\}$, $\{(b,ab),(aba,ba),(ba,a),(ab,aba)\}$, $\{(ab,b),(b,b),(aba,b),(ba,b)\}$, $\{(ab,ab),(b,ba),(aba,a),(ba,aba)\}$, $\{(ba,ba),(ab,a),(b,aba),(aba,ab)\}$, and $\{(aba,aba),(ba,ab),(ab,ba),(b,a)\}$. If $C$ is any one of these six 4-cycles, then no two distinct elements $x$ and $y$ of $C$ generates Q1×Q1, because $\{x,y\} \subseteq C$ and $C$ is a proper subquadratical quasigroup of Q1×Q1, isomorphic to Q1.

**Example 4.5.** Using the direct product, we can construct quadratical quasigroups of various orders such as $5^n$, $9^n$, $(45)^n$ and so on. Using Dudek's quadratical quasigroups of orders 13 and 17 we can construct even more quadratical quasigroups of orders greater than those considered in [1].

**Example 4.6.** $(Q1 \times Q1)^* = (Q1)^* \times (Q1)^*$ and $[Q1 \times (Q1)^*]^* = (Q1)^* \times Q1$. Note that $(a,ba)$ and $(ab,b)$ generate $Q1 \times (Q1)^*$ and $(ba,a)$ and $(b,ab)$ generate $(Q1)^* \times Q1$ while $Q1 \times Q1$ and $(Q1)^* \times (Q1)^*$ are **not** 2-generated.

## 5. The elements $nk^*$.

The following Theorem is easily proved for $k = 1$ and, by induction on $k$, is straightforward to prove for all $k \in \{0,1,2,3,...\} = \mathbb{N}_0$. The proof is omitted.

**Theorem 5.1.** *For all* $k \in \mathbb{N}_0$,
$[(1+4k)1]^* = (1+4k)1$; $[(1+4k)2]^* = (1+4k)3$; $[(1+4k)3]^* = (1+4k)2$; $[(1+4k)4]^* = (1+4k)4$;
$[(2+4k)1]^* = (2+4k)3$; $[(2+4k)2]^* = (2+4k)4$; $[(2+4k)3]^* = (2+4k)1$; $[(2+4k)4]^* = (2+4k)2$;
$[(3+4k)1]^* = (3+4k)4$; $[(3+4k)2]^* = (3+4k)2$; $[(3+4k)3]^* = (3+4k)3$; $[(3+4k)4]^* = (3+4k)1$ *and*
$[(4+4k)1]^* = (4+4k)2$; $[(4+4k)2]^* = (4+4k)1$; $[(4+4k)3]^* = (4+4k)4$; $[(4+4k)4]^* = (4+4k)3$.

Now, considering the quadratical quasigroups of form $Qn$, from the remarks in the paragraph preceding Example 4.1, we see that there are at most 4 groupoids of the form $Qn$ for any given integer $n$. Since the dual of a quadratical quasigroup of the form $Qn$ must also have the form $Qn$, we can tell, from the following Theorem, which values of $aba \cdot a$ may yield groupoids that are duals of each other.

**Theorem 5.2.** *For all positive integers* $n \geq 2$, *the following identities are valid in a quadratical quasigroup of form* $Qn$, *depending on the value of* $aba \cdot a$:

| $aba \cdot a$ | $aba \cdot ab$ | $aba \cdot ba$ | $aba \cdot b$ | $a \cdot aba$ | $ab \cdot aba$ | $ba \cdot aba$ | $b \cdot aba$ | $n1 \cdot n2$ | $n2 \cdot n4$ | $n3 \cdot n1$ | $n4 \cdot n3$ |
|---|---|---|---|---|---|---|---|---|---|---|---|
| $n1$ | $n2$ | $n3$ | $n4$ | $n2$ | $n4$ | $n1$ | $n3$ | $a$ | $ab$ | $ba$ | $b$ |
| $n2$ | $n4$ | $n1$ | $n3$ | $n4$ | $n3$ | $n2$ | $n1$ | $ba$ | $a$ | $b$ | $ab$ |
| $n3$ | $n1$ | $n4$ | $n2$ | $n1$ | $n2$ | $n3$ | $n4$ | $ab$ | $b$ | $a$ | $ba$ |
| $n4$ | $n3$ | $n2$ | $n1$ | $n3$ | $n1$ | $n4$ | $n2$ | $b$ | $ba$ | $ab$ | $a$ |

| $aba \cdot a$ | $11 \cdot 34$ | $23 \cdot 14$ | $34 \cdot 14$ | $14 \cdot 21$ | |
|---|---|---|---|---|---|
| $n1$ | $n3$ | $n2$ | $n1$ | $n1$ | $(n-1)2 = 11 \cdot n1 = n2 \cdot 11$ |
| $n2$ | $n1$ | $n4$ | $n2$ | $n2$ | $(n-1)4 = 11 \cdot n2 = n4 \cdot 11$ |
| $n3$ | $n4$ | $n1$ | $n3$ | $n3$ | $(n-1)1 = 11 \cdot n3 = n1 \cdot 11$ |
| $n4$ | $n2$ | $n3$ | $n4$ | $n4$ | $(n-1)3 = 11 \cdot n4 = n3 \cdot 11$ |

Proof of Theorem 5.2: We prove only the identities for when $aba \cdot a = n2$, as the proofs of the other three cases are similar. We have $aba \cdot a = n2$ Then, $aba = (a \cdot aba)(aba \cdot a) = (a \cdot aba)(n2)$. By Proposition 3.11 and Theorem 2.1, $a \cdot aba = n4 = a \cdot bab = aba \cdot ab$. Then, $n4 = aba \cdot ab = (aba \cdot a)(aba \cdot b) = n2 \cdot (aba \cdot b)$. By Proposition 3.5, $aba \cdot b = n3$. So $aba \cdot ba = (aba \cdot b)(aba \cdot a) = n3 \cdot n2 = n1$ (by Proposition 3.5). Then, $aba = (b \cdot aba)(aba \cdot b) = (b \cdot aba) \cdot n3$, which by Proposition 3.11 implies $b \cdot aba = n1$. Then, using 3.5, $ab \cdot aba = (a \cdot aba)(b \cdot aba) = n4 \cdot n1 = n3$ and $ba \cdot aba = (b \cdot aba)(a \cdot aba) = n1 \cdot n4 = n2$. We also have $n1 \cdot n2 = (aba \cdot ba)(ba \cdot aba) = ba$, $n2 \cdot n4 = (aba \cdot a)(a \cdot aba) = a$, $n3 \cdot n1 = (aba \cdot b)(b \cdot aba) = b$ and $n4 \cdot n3 = (aba \cdot ab)(ab \cdot aba) = ab$. Now,

$11 \cdot 34 = a \cdot [(b \cdot ba)(ba \cdot a)] = [a \cdot (b \cdot ba)][a \cdot (ba \cdot a)] = (ab \cdot aba)(aba \cdot a) = n3 \cdot n2 = n1$,

$34 \cdot 14 = [(b \cdot ba)(ba \cdot a)] \cdot b = [(b \cdot ba) \cdot b][(ba \cdot a) \cdot b] = (b \cdot bab)(bab \cdot ab) = (b \cdot aba)(aba \cdot ab) = n1 \cdot n4 = n2$,

$14 \cdot 21 = b \cdot (a \cdot ab) = ba \cdot bab = ba \cdot aba = n2$ and $23 \cdot 14 = (ba \cdot a) \cdot b = bab \cdot ab = aba \cdot ab = n4$.

Finally, $a \cdot n2 = a \cdot aba \cdot a = aba \cdot a \cdot aba = aba \cdot n4 = (n-1)4 = 11 \cdot n2$ and $n4 \cdot a = a \cdot aba \cdot a = aba \cdot a \cdot aba = (n-1)4 = n4 \cdot 11$. This completes the proof of the validity of the identities indicated in row 3 of the two tables in Theorem 5.2, when $aba \cdot a = n2$. ∎

As mentioned above, Theorem 5.2 will be useful when we look for the duals of the quadratical quasigroups that we will call Q3 and Q4, as will the following concept.

**Definition 5.1.** If a quadratical quasigroup of form *Qn* exists for some integer *n* then ***the identity generated on the left*** [***on the right***] ***by an identity*** $kr \cdot ls = mt$ $(k,l,m \le n, \{r,s,t\} \subseteq \{1,2,3,4\})$ is defined as the identity $(aba \cdot kr) \cdot (aba \cdot ls) = aba \cdot mt$ [$(kr \cdot aba) \cdot (ls \cdot aba) = mt \cdot aba$] and $kr \cdot ls = mt$ is called the ***generating identity.***

Note that Propositions 3.6 and 3.7, along with Theorem 4.2, give the means of calculating identities generated on the left and right by a given identity. Multiplying on the left (or on the right) repeatedly *n*-times gives *n* distinct identities. These methods will later be used to prove that quadratical quasigroups of the form *Q*6 do not exist.

## 6. Quadratical quasigroups of forms *Q*3 and *Q*4.

We give the Cayley tables of Dudek's quadratical quasigroups of orders 13 and 17. He proved that there were two quadratical quasigroups of each order and we prove that each pair of groupoids are duals of each other.

First we note that for a quadratical quasigroup of form $Q3$, if $aba \cdot a = n3 = 33 = (ba \cdot a)(a \cdot ab)$ then $aba \cdot a = a(ba \cdot a) = (a \cdot ab)(aba) = ab$ which implies, by cancellation, $ba \cdot a = b$, a contradiction because $H1 \cap H2 = \emptyset$. If $aba \cdot a = n4 = 34 = (b \cdot ba)(ba \cdot a)$ then $ab \cdot aba = a \cdot (b \cdot ba) = (ba \cdot a)(aba) = a$, which implies $b \cdot ba = a$, a contradiction. Hence, $aba \cdot a \in \{31, 32\} = \{a \cdot ab, ab \cdot b\}$. Setting $aba \cdot a = a \cdot ab$ and using the properties of quadratical quasigroups (Theorem 2.1) we obtain the Cayley Table 3. It can be checked that it is medial and bookend and so, by Theorem 2.11, this groupoid that we call Q3 is a quadratical quasigroup.

**Table 3.**

| Q3 | 11 | 12 | 13 | 14 | aba | 21 | 22 | 23 | 24 | 31 | 32 | 33 | 34 |
|---|---|---|---|---|---|---|---|---|---|---|---|---|---|
| 11 | **11** | **21** | *aba* | **12** | **32** | **14** | **23** | **31** | **34** | **22** | **13** | **24** | **33** |
| 12 | *aba* | **12** | **14** | **22** | **34** | **32** | **13** | **33** | **21** | **24** | **23** | **31** | **11** |
| 13 | **23** | **11** | **13** | *aba* | **31** | **24** | **32** | **12** | **33** | **14** | **34** | **21** | **22** |
| 14 | **13** | *aba* | **24** | **14** | **33** | **31** | **34** | **22** | **11** | **32** | **21** | **12** | **23** |
| *aba* | **31** | **32** | **33** | **34** | *aba* | **11** | **12** | **13** | **14** | **21** | **22** | **23** | **24** |
| 21 | **32** | **23** | **34** | **13** | **22** | **21** | **31** | *aba* | **22** | **24** | **33** | **11** | **14** |
| 22 | **33** | **34** | **11** | **21** | **24** | *aba* | **22** | **24** | **32** | **12** | **23** | **13** | **31** |
| 23 | **24** | **14** | **31** | **32** | **21** | **33** | **21** | **23** | *aba* | **34** | **12** | **22** | **13** |
| 24 | **12** | **31** | **22** | **33** | **23** | **23** | *aba* | **34** | **24** | **11** | **14** | **32** | **21** |
| 31 | **34** | **13** | **21** | **24** | **22** | **12** | **33** | **14** | **23** | **31** | **11** | *aba* | **32** |
| 32 | **22** | **33** | **23** | **11** | **24** | **13** | **14** | **21** | **31** | *aba* | **32** | **34** | **12** |
| 33 | **14** | **22** | **32** | **23** | **21** | **34** | **24** | **11** | **12** | **13** | **31** | **33** | *aba* |
| 34 | **21** | **24** | **12** | **31** | **23** | **22** | **11** | **32** | **13** | **33** | *aba* | **14** | **34** |

There are then two ways to obtain the Cayley table for $(Q3)^*$. Firstly, we can use $aba * a = 32^* = (a*b)*b$ and, using the properties of quadratical quasigroups, we can then calculate the remaining products in Table 4.

Alternatively, we can calculate the products directly from Table 3, using Theorem 5.1. For example, $(23)^* = (b*a)*a = a(ab) = 21$, and $(32)^* = [(a*b)*b)] * [b*(b*a)] = (ab \cdot b) \cdot (b \cdot ba) = 32$. Hence, $(32)^* * (23)^* = 32 * 21 = 21 \cdot 32$. From Table 3, $21 \cdot 32 = 33$. But from Theorem 5.1, $33 = 33^*$. So we obtain $(32)^* * (23)^* = 33 = 33^*$. The remaining products in Table 4 can be calculated in similar fashion. Having already checked that Table 3 is quadratical, Table 4 also produces a quadratical quasigroup, the dual groupoid $(Q3)^*$.

**Table 4.**

| (Q3)* | 11* | 12* | 13* | 14* | aba | 21* | 22* | 23* | 24* | 31* | 32* | 33* | 34* |
|---|---|---|---|---|---|---|---|---|---|---|---|---|---|
| 11* | **11*** | **21*** | *aba* | **12*** | **34*** | **22*** | **13*** | **32*** | **33*** | **23*** | **24*** | **14*** | **31*** |
| 12* | *aba* | **12*** | **14*** | **22*** | **33*** | **34*** | **24*** | **31*** | **14*** | **13*** | **21*** | **32*** | **23*** |
| 13* | **23*** | **11*** | **13*** | *aba* | **32*** | **14*** | **34*** | **21*** | **31*** | **22*** | **33*** | **24*** | **12*** |
| 14 * | **13*** | *aba* | **24*** | **14*** | **31*** | **32*** | **33*** | **12*** | **23*** | **34*** | **14*** | **21*** | **22*** |
| aba | **32*** | **34*** | **31*** | **33*** | *aba* | **11*** | **12*** | **13*** | **14*** | **21*** | **22*** | **23*** | **24*** |
| 21* | **34*** | **23*** | **33*** | **24*** | **12*** | **21*** | **31*** | *aba* | **22*** | **23*** | **23*** | **11*** | **14*** |
| 22* | **31*** | **33*** | **23*** | **14*** | **14*** | *aba* | **22*** | **24*** | **32*** | **34*** | **34*** | **13*** | **21*** |
| 23* | **14*** | **22*** | **32*** | **34*** | **11*** | **33*** | **21*** | **23*** | *aba* | **12*** | **12*** | **31*** | **13*** |
| 24* | **21*** | **32*** | **12*** | **31*** | **13*** | **23*** | *aba* | **34*** | **24*** | **14*** | **14*** | **22*** | **33*** |
| 31* | **33*** | **24*** | **14*** | **21*** | **22*** | **12*** | **23*** | **14*** | **34*** | **31*** | **13*** | *aba* | **32*** |
| 32* | **12*** | **31*** | **22*** | **23*** | **24*** | **13*** | **14*** | **33*** | **21*** | *aba* | **32*** | **34*** | **11*** |
| 33* | **22*** | **23*** | **34*** | **13*** | **21*** | **24*** | **32*** | **14*** | **12*** | **14*** | **31*** | **33*** | *aba* |
| 34* | **24*** | **14*** | **21*** | **32*** | **23*** | **31*** | **11*** | **22*** | **13*** | **33*** | *aba* | **12*** | **34*** |

Similarly, we can calculate the Cayley Tables for Q4 and its dual $(Q4)^*$:

**Table 5.**

| Q4 | 11 | 12 | 13 | 14 | aba | 21 | 22 | 23 | 24 | 31 | 32 | 33 | 34 | 41 | 42 | 43 | 44 |
|---|---|---|---|---|---|---|---|---|---|---|---|---|---|---|---|---|---|
| 11 | **11** | **21** | *aba* | **12** | **44** | **24** | **32** | **42** | **43** | **14** | **23** | **31** | **41** | **33** | **34** | **13** | **22** |
| 12 | *aba* | **12** | **14** | **22** | **43** | **44** | **23** | **41** | **34** | **32** | **13** | **42** | **21** | **11** | **31** | **24** | **33** |
| 13 | **23** | **11** | **13** | *aba* | **42** | **31** | **44** | **22** | **41** | **24** | **43** | **12** | **33** | **32** | **21** | **34** | **14** |
| 14 | **13** | *aba* | **24** | **14** | **41** | **42** | **43** | **33** | **21** | **44** | **34** | **22** | **11** | **23** | **12** | **31** | **32** |
| aba | **42** | **44** | **41** | **43** | *aba* | **11** | **12** | **13** | **14** | **21** | **22** | **23** | **24** | **31** | **32** | **33** | **34** |
| 21 | **44** | **32** | **43** | **23** | **12** | **21** | **31** | *aba* | **22** | **34** | **42** | **11** | **14** | **24** | **33** | **41** | **13** |
| 22 | **41** | **43** | **21** | **34** | **14** | *aba* | **22** | **24** | **32** | **12** | **33** | **13** | **44** | **42** | **23** | **11** | **31** |
| 23 | **31** | **24** | **42** | **44** | **11** | **33** | **21** | **23** | *aba* | **41** | **12** | **32** | **13** | **34** | **14** | **22** | **43** |
| 24 | **22** | **42** | **33** | **41** | **13** | **23** | *aba* | **34** | **24** | **11** | **14** | **43** | **31** | **12** | **44** | **32** | **21** |
| 31 | **43** | **23** | **34** | **13** | **22** | **12** | **42** | **14** | **33** | **31** | **41** | *aba* | **32** | **44** | **11** | **21** | **24** |
| 32 | **33** | **41** | **11** | **21** | **24** | **13** | **14** | **31** | **44** | *aba* | **32** | **34** | **42** | **22** | **43** | **23** | **12** |
| 33 | **24** | **14** | **44** | **32** | **21** | **41** | **34** | **11** | **12** | **43** | **31** | **33** | *aba* | **13** | **22** | **42** | **23** |
| 34 | **12** | **31** | **22** | **42** | **23** | **32** | **11** | **43** | **13** | **33** | *aba* | **44** | **34** | **21** | **24** | **14** | **41** |
| 41 | **21** | **34** | **12** | **31** | **32** | **14** | **33** | **44** | **23** | **22** | **11** | **24** | **43** | **41** | **13** | *aba* | **42** |
| 42 | **14** | **22** | **32** | **33** | **34** | **43** | **13** | **21** | **31** | **23** | **24** | **41** | **12** | *aba* | **42** | **44** | **11** |
| 43 | **32** | **33** | **23** | **11** | **31** | **34** | **24** | **12** | **42** | **13** | **44** | **21** | **22** | **14** | **41** | **43** | *aba* |
| 44 | **34** | **13** | **31** | **24** | **33** | **22** | **41** | **32** | **11** | **42** | **21** | **14** | **23** | **43** | *aba* | **12** | **44** |

**Table 6.**

| (Q4)* | 11* | 12* | 13* | 14* | aba | 21* | 22* | 23* | 24* | 31* | 32* | 33* | 34* | 41* | 42* | 43* | 44* |
|---|---|---|---|---|---|---|---|---|---|---|---|---|---|---|---|---|---|
| 11* | 11* | 21* | aba | 12* | 41* | 34* | 24* | 43* | 42* | 13* | 33* | 22* | 44* | 14* | 23* | 31* | 32* |
| 12* | aba | 12* | 14* | 22* | 42* | 41* | 33* | 44* | 23* | 24* | 11* | 43* | 31* | 32* | 13* | 34* | 21* |
| 13* | 23* | 11* | 13* | aba | 43* | 22* | 41* | 32* | 44* | 34* | 42* | 14* | 21* | 24* | 31* | 12* | 33* |
| 14* | 13* | aba | 24* | 14* | 44* | 43* | 42* | 21* | 31* | 41* | 23* | 32* | 12* | 33* | 34* | 22* | 11* |
| aba | 43* | 41* | 44* | 42* | aba | 11* | 12* | 13* | 14* | 21* | 22* | 23* | 24* | 31* | 32* | 33* | 34* |
| 21* | 41* | 24* | 42* | 33* | 12* | 21* | 31* | aba | 22* | 44* | 34* | 11* | 14* | 23* | 43* | 32* | 13* |
| 22* | 44* | 42* | 31* | 23* | 14* | aba | 22* | 24* | 32* | 12* | 43* | 13* | 33* | 34* | 21* | 11* | 41* |
| 23* | 22* | 34* | 43* | 41* | 11* | 33* | 21* | 23* | aba | 32* | 12* | 42* | 13* | 44* | 14* | 24* | 31* |
| 24* | 32* | 43* | 21* | 44* | 13* | 23* | aba | 34* | 24* | 11* | 14* | 31* | 41* | 12* | 33* | 42* | 22* |
| 31* | 42* | 33* | 23* | 11* | 22* | 12* | 34* | 14* | 43* | 31* | 41* | aba | 32* | 13* | 44* | 21* | 24* |
| 32* | 21* | 44* | 12* | 31* | 24* | 13* | 14* | 41* | 33* | aba | 32* | 34* | 42* | 22* | 11* | 23* | 43* |
| 33* | 34* | 13* | 41* | 24* | 21* | 32* | 44* | 11* | 12* | 43* | 31* | 33* | aba | 42* | 22* | 14* | 23* |
| 34* | 14* | 22* | 32* | 43* | 23* | 42* | 11* | 31* | 13* | 33* | aba | 44* | 34* | 21* | 24* | 41* | 12* |
| 41* | 31* | 23* | 34* | 13* | 32* | 14* | 43* | 33* | 21* | 22* | 44* | 24* | 11* | 41* | 12* | aba | 42* |
| 42* | 33* | 32* | 11* | 21* | 34* | 31* | 13* | 22* | 41* | 23* | 24* | 12* | 43* | aba | 42* | 44* | 14* |
| 43* | 24* | 14* | 33* | 32* | 31* | 44* | 23* | 12* | 34* | 42* | 13* | 21* | 22* | 11* | 41* | 43* | aba |
| 44* | 12* | 31* | 22* | 34* | 33* | 24* | 32* | 42* | 11* | 14* | 21* | 41* | 23* | 43* | aba | 13* | 44* |

## 7. No quadratical quasigroup of form *Q*6 exists

**Theorem 7.1.** *There is no quadratical quasigroup of form Q6.*

Proof: Case 1: $aba \cdot a = 61$. Using 3.6, 3.7, 5.2, alterability and idempotency we have

(1): $aba \cdot a = 61 \Rightarrow a \cdot 61 = 61 \cdot aba = 52 = 62 \cdot a = aba \cdot 62 = 52 \cdot 52$.

By 3.5 and 3.6, $62 = 61 \cdot 64 \Rightarrow 52 = 51 \cdot 54$ and by Definition 3.1, $62 = 52 \cdot 54 \Rightarrow 52 = 42 \cdot 44$ and so

(2): $52 = 51 \cdot 54 = 42 \cdot 44$.

By Theorem 5.2, $61 = 14 \cdot 21 = 34 \cdot 14$, $62 = 23 \cdot 14$ and $63 = 11 \cdot 34$ and these identities generate the following:

(3): $52 = 63 \cdot 12 = 13 \cdot 64 = 64 \cdot 21 = 23 \cdot 13$.

Using (1), (2) and (3), we see that for the identity $52 = 12 \cdot x$, $x \in \{14, 22, 23, 24, 31, 32, 33, 34, 41, 42, 43, 51, 53, 63\}$. Now by definition 3.1, $22 = 12 \cdot 14 \neq 52$ and so $x \neq 14$.

To eliminate the other possibilities for *x* we now use the generating identities (4) through (14), indicated in the chart below. Assuming that *Q*6 is quadratical, using the properties of a quadratical quasigroup we will prove that all the remaining possible values of *x* lead to a contradiction.

**Chart 1.**

| | (4) $(n-1)2$ $=11 \cdot n1$ | (5) $(n-1)2$ $=n2 \cdot 11$ | (6) $aba \cdot 11$ $= 61$ | (7) $11 \cdot aba$ $= 62$ | (8) Propn. 3.5 | (9) Defn. 3.1 | (10) $n1 =$ $14 \cdot 21$ | (11) $n2 =$ $23 \cdot 14$ | (12) $n3 =$ $11 \cdot 34$ | (13) $n1 =$ $34 \cdot 14$ | (14) *idem* |
|---|---|---|---|---|---|---|---|---|---|---|---|
| 52 | $11 \cdot 61$ | $62 \cdot 11$ | $aba \cdot 62$ | $61 \cdot aba$ | $51 \cdot 54$ | $42 \cdot 44$ | $63 \cdot 12$ | $13 \cdot 64$ | $64 \cdot 21$ | $23 \cdot 13$ | $52 \cdot 52$ |
| 44 | $62 \cdot 52$ | $54 \cdot 62$ | $aba \cdot 54$ | $52 \cdot aba$ | $42 \cdot 43$ | $34 \cdot 33$ | $51 \cdot 64$ | $61 \cdot 22$ | $53 \cdot 12$ | $11 \cdot 61$ | $44 \cdot 44$ |
| 33 | $54 \cdot 44$ | $43 \cdot 54$ | $aba \cdot 43$ | $44 \cdot aba$ | $34 \cdot 31$ | $23 \cdot 21$ | $42 \cdot 53$ | $52 \cdot 41$ | $41 \cdot 64$ | $62 \cdot 52$ | $33 \cdot 33$ |
| 21 | $43 \cdot 33$ | $31 \cdot 43$ | $aba \cdot 31$ | $33 \cdot aba$ | $23 \cdot 22$ | $11 \cdot 12$ | $34 \cdot 41$ | $44 \cdot 32$ | $32 \cdot 53$ | $54 \cdot 34$ | $21 \cdot 21$ |
| 12 | $31 \cdot 21$ | $22 \cdot 31$ | $aba \cdot 22$ | $21 \cdot aba$ | $11 \cdot 14$ | $62 \cdot 64$ | $23 \cdot 32$ | $33 \cdot 24$ | $24 \cdot 41$ | $43 \cdot 33$ | $12 \cdot 12$ |
| 64 | $22 \cdot 12$ | $14 \cdot 22$ | $aba \cdot 14$ | $12 \cdot aba$ | $62 \cdot 63$ | $54 \cdot 53$ | $11 \cdot 24$ | $21 \cdot 13$ | $13 \cdot 32$ | $31 \cdot 21$ | $64 \cdot 64$ |
| 53 | $14 \cdot 64$ | $63 \cdot 14$ | $aba \cdot 63$ | $64 \cdot aba$ | $54 \cdot 51$ | $43 \cdot 41$ | $62 \cdot 13$ | $12 \cdot 61$ | $61 \cdot 24$ | $22 \cdot 12$ | $53 \cdot 53$ |
| 41 | $63 \cdot 53$ | $51 \cdot 63$ | $aba \cdot 51$ | $53 \cdot aba$ | $43 \cdot 42$ | $31 \cdot 32$ | $54 \cdot 61$ | $64 \cdot 52$ | $52 \cdot 13$ | $14 \cdot 54$ | $41 \cdot 41$ |
| 32 | $51 \cdot 41$ | $42 \cdot 51$ | $aba \cdot 42$ | $41 \cdot aba$ | $31 \cdot 34$ | $22 \cdot 24$ | $43 \cdot 52$ | $53 \cdot 44$ | $44 \cdot 61$ | $63 \cdot 53$ | $32 \cdot 32$ |
| 24 | $42 \cdot 32$ | $34 \cdot 42$ | $aba \cdot 34$ | $32 \cdot aba$ | $22 \cdot 23$ | $14 \cdot 13$ | $31 \cdot 44$ | $41 \cdot 33$ | $33 \cdot 52$ | $51 \cdot 41$ | $24 \cdot 24$ |
| 13 | $34 \cdot 24$ | $23 \cdot 34$ | $aba \cdot 23$ | $24 \cdot aba$ | $14 \cdot 11$ | $63 \cdot 61$ | $22 \cdot 33$ | $32 \cdot 21$ | $21 \cdot 44$ | $42 \cdot 32$ | $13 \cdot 13$ |
| 61 | $23 \cdot 13$ | $11 \cdot 23$ | $aba \cdot 11$ | $13 \cdot aba$ | $63 \cdot 62$ | $51 \cdot 52$ | $14 \cdot 21$ | $24 \cdot 12$ | $12 \cdot 33$ | $34 \cdot 14$ | $61 \cdot 61$ |
| 51 | $13 \cdot 63$ | $61 \cdot 13$ | $aba \cdot 61$ | $63 \cdot aba$ | $53 \cdot 52$ | $41 \cdot 42$ | $64 \cdot 11$ | $14 \cdot 62$ | $62 \cdot 23$ | $24 \cdot 64$ | $51 \cdot 51$ |
| 31 | $53 \cdot 43$ | $41 \cdot 53$ | $aba \cdot 41$ | $43 \cdot aba$ | $33 \cdot 32$ | $21 \cdot 22$ | $44 \cdot 51$ | $54 \cdot 42$ | $42 \cdot 63$ | $64 \cdot 44$ | $31 \cdot 31$ |
| 11 | $33 \cdot 23$ | $21 \cdot 33$ | $aba \cdot 21$ | $23 \cdot aba$ | $13 \cdot 12$ | $61 \cdot 62$ | $24 \cdot 31$ | $34 \cdot 22$ | $22 \cdot 43$ | $44 \cdot 24$ | $11 \cdot 11$ |
| 62 | $21 \cdot 11$ | $12 \cdot 21$ | $aba \cdot 12$ | $11 \cdot aba$ | $61 \cdot 64$ | $52 \cdot 54$ | $13 \cdot 22$ | $23 \cdot 14$ | $14 \cdot 31$ | $33 \cdot 23$ | $62 \cdot 62$ |
| 54 | $12 \cdot 62$ | $64 \cdot 12$ | $aba \cdot 64$ | $62 \cdot aba$ | $52 \cdot 53$ | $44 \cdot 43$ | $61 \cdot 14$ | $11 \cdot 63$ | $63 \cdot 22$ | $21 \cdot 11$ | $54 \cdot 54$ |
| 43 | $64 \cdot 54$ | $53 \cdot 64$ | $aba \cdot 53$ | $54 \cdot aba$ | $44 \cdot 41$ | $33 \cdot 31$ | $52 \cdot 63$ | $62 \cdot 51$ | $51 \cdot 14$ | $12 \cdot 62$ | $43 \cdot 43$ |
| 34 | $52 \cdot 42$ | $44 \cdot 52$ | $aba \cdot 44$ | $42 \cdot aba$ | $32 \cdot 33$ | $24 \cdot 23$ | $41 \cdot 54$ | $51 \cdot 43$ | $43 \cdot 62$ | $61 \cdot 51$ | $34 \cdot 34$ |
| 23 | $44 \cdot 34$ | $33 \cdot 44$ | $aba \cdot 33$ | $34 \cdot aba$ | $24 \cdot 21$ | $13 \cdot 11$ | $32 \cdot 43$ | $42 \cdot 31$ | $31 \cdot 54$ | $52 \cdot 42$ | $23 \cdot 23$ |
| 14 | $32 \cdot 22$ | $24 \cdot 32$ | $aba \cdot 24$ | $22 \cdot aba$ | $12 \cdot 13$ | $64 \cdot 63$ | $21 \cdot 34$ | $31 \cdot 23$ | $23 \cdot 42$ | $41 \cdot 31$ | $14 \cdot 14$ |
| 63 | $24 \cdot 14$ | $13 \cdot 24$ | $aba \cdot 13$ | $14 \cdot aba$ | $64 \cdot 61$ | $53 \cdot 51$ | $12 \cdot 23$ | $22 \cdot 11$ | $11 \cdot 34$ | $32 \cdot 22$ | $63 \cdot 63$ |
| 42 | $61 \cdot 51$ | $52 \cdot 61$ | $aba \cdot 52$ | $51 \cdot aba$ | $41 \cdot 44$ | $32 \cdot 34$ | $53 \cdot 62$ | $63 \cdot 54$ | $54 \cdot 11$ | $13 \cdot 63$ | $42 \cdot 42$ |
| 22 | $41 \cdot 31$ | $32 \cdot 41$ | $aba \cdot 32$ | $31 \cdot aba$ | $21 \cdot 24$ | $12 \cdot 14$ | $33 \cdot 42$ | $43 \cdot 34$ | $34 \cdot 51$ | $53 \cdot 53$ | $22 \cdot 22$ |

When we use a particular value of an element we will refer to the column in which this value appears, such as (4), (5), …, (13) or (14). For example, if $52 = 12 \cdot 63$ then, using (10), $52 = 12 \cdot 63 = 63 \cdot 12$, a contradiction. We will use the fact that $52 = 63 \cdot 12$, from (10), henceforth without mention.

By (10), if $52 = 12 \cdot 53 = 63 \cdot 12$ then $12 = 53 \cdot 63$ and, multiplying on the right by *aba* gives $64 = 41 \cdot 51$ which, along with $51 \cdot 41 = 24$ (from (13)) gives $51 = 64 \cdot 24$. This contradicts $51 = 64 \cdot 11$, from (10).

If $52 = 12 \cdot 51 = 63 \cdot 12$ then $12 = 51 \cdot 63 = 62 \cdot 64$, from (9). Hence, by (8) and (9), $61 = 63 \cdot 62 = 64 \cdot 51 = 51 \cdot 52$. Therefore, using (13), $51 = 52 \cdot 64 = 24 \cdot 64$, a contradiction.

If $52 = 12 \cdot 43 = 63 \cdot 12$ then, by (12), $12 = 43 \cdot 63 = 24 \cdot 41$. By 3.11, $63 \cdot 24 = 41 \cdot 43 = aba = 23 \cdot 24$, contradiction.

If $52 = 12 \cdot 42 = 63 \cdot 12$ then, by (12), $12 = 42 \cdot 63 = 24 \cdot 41$. By 3.1 and (13), $51 = 41 \cdot 42 = 63 \cdot 24 = 24 \cdot 64$. So, by (9), $24 = 64 \cdot 63 = 14$, contradiction.

If $52 = 12 \cdot 41 = 63 \cdot 12$ then, by (12), $12 = 41 \cdot 63 = 24 \cdot 41$ and so, using (4), $41 = 63 \cdot 24 = 63 \cdot 53$, contradiction.

If $52 = 12 \cdot 34 = 63 \cdot 12$ then, by (10), $12 = 34 \cdot 63 = 23 \cdot 32$ and, by 3.1 and (11), $42 = 32 \cdot 34 = 63 \cdot 23 = 63 \cdot 54$, contradiction.

If $52 = 12 \cdot 33 = 63 \cdot 12$ then, by (10), $12 = 33 \cdot 63 = 23 \cdot 32$ and so, by 3.1 and 3.5, $34 = 32 \cdot 33 = 63 \cdot 23 = 24 \cdot 23$, contradiction.

If $52 = 12 \cdot 32 = 63 \cdot 12$ then, by (10), $12 = 32 \cdot 63 = 23 \cdot 32$ and so, by (13), $32 = 63 \cdot 23 = 63 \cdot 53$, contradiction.

If $52 = 12 \cdot 31 = 63 \cdot 12$ then, by (4), $12 = 31 \cdot 63 = 31 \cdot 21$, contradiction.

If $52 = 12 \cdot 24 = 63 \cdot 12$ then, by (4), $12 = 24 \cdot 63 = 24 \cdot 41$, contradiction.

If $52 = 12 \cdot 23 = 63 \cdot 12$ then, by (10), $12 = 23 \cdot 63 = 23 \cdot 32$, contradiction.

If $52 = 12 \cdot 22 = 63 \cdot 12$ then, by (5), $12 = 22 \cdot 63 = 22 \cdot 31$, contradiction.

If $52 = 12 \cdot 14$ then, by 3.1, $52 = 12 \cdot 14 = 22$, contradiction.

We have proved that when $aba \cdot a = 61$, there is no right solvability, a contradiction. The proof that there is no right solvability when $aba \cdot a \in \{62, 63, 64\}$ is similar, where the values in Chart 1 are different, according to Theorem 5.2. We omit these detailed calculations. ∎

So far we have seen that there are quadratical quasigroups of orders 1,5,9,13,17,25 and any finite product of these numbers. Suppose $A = \{1, 5, 9, 13, 17, 25\}$ and $B = \{x : x \text{ is a finite product of members of } A\}$.

**Definition 7.1.** The spectrum *SpV* of a variety *V* is the set of all orders of finite members of *V*.

From Proposition 3.4 and the fact that there is no quadratical quasigroup of order 21 [1], we have:

**Corollary 7.2.** $B \subseteq Sp\ Q \subseteq \{1 + 4n : n \in \{0, 1, 2, 3, 4, 6, 7, 8, 9, ...\}\}$

## 8. Translatable groupoids

Patterns of "translatability" can be hidden in the Cayley tables of quadratical quasigroups. One can assume the properties of quadratical quasigroups and then calculate whether translatable groupoids of various orders exist with these properties. We proceed to prove that the quadratical quasigroups Q1, (Q1)*, Q3, (Q3)*, Q4 and (Q4)* are translatable and that Q2 is *not* translatable.

**Definition 8.1.** A finite groupoid $G = \{1, 2, 3, ..., n-1, n\}$ is called *k-translatable* $(k \in \{1, 2, ..., n\})$ if its Cayley Table is obtained by the following rule: If the first row of the Cayley Table is:

|   | 1 | 2 | 3 | . | . | . | n-2 | n-1 | n |
|---|---|---|---|---|---|---|-----|-----|---|
| 1 | $x_1$ | $x_2$ | $x_3$ | . | . | . | $x_{n-2}$ | $x_{n-1}$ | $x_n$ |

then the q-th row is obtained from the (q-1)-st row by taking the last k entries in the (q-1)-st row and inserting them as the first k entries of the q-th row and by taking the first *n*-k entries of the (q-1)-st row and inserting them as the last (*n*-k) entries of the q-th row, where $q \in \{2,3,4,...,n\}$. Then the (ordered) sequence $x_1, x_2, x_3,..., x_n$ is called a **k-translatable sequence of G** = **{1, 2, ... ,n}**. A groupoid is called **a translatable groupoid** if it has a k-translatable sequence for some $k \in \{1,2,...,n\}$.

Note that in this definition we adopt the usual convention that the left hand column of the table, from top to bottom, is the top row from left to right. That is, a Cayley table for a groupoid $G = \{1,2,3,...,n\}$ is a product table:

| G | 1 | 2 | 3 | . | . | . | $n-2$ | $n-1$ | n |
|---|---|---|---|---|---|---|---|---|---|
| 1 | | | | | | | | | |
| 2 | | | | | | | | | |
| 3 | | | | | | | | | |
| . | | | | | | | | | |
| . | | | | | | | | | |
| . | | | | | | | | | |
| $n-2$ | | | | | | | | | |
| $n-1$ | | | | | | | | | |
| n | | | | | | | | | |

It is important to note that a k-translatable sequence of a groupoid depends on the ordering 1, 2, ... , n of the elements of *G*. A groupoid may be k-translatable for one ordering but not for another (see example 8.4 below).

**Example 8.1.** The groupoid $\{1,2,3,...,n\}$ with product as addition mod-*n* is (*n*-1)-translatable.

**Proposition 8.1.** *No idempotent groupoid* G *of order greater than 1 is 1-translatable.*

Proof: If $G = \{1,2,3,...,n\}$ is idempotent and 1-translatable then by definition, $2 \cdot 2 = 1$, a contradiction. ∎

**Example 8.2.** The following idempotent groupoid of order *n* = 5 is (*n*-1)-translatable.

|   | 1 | 2 | 3 | 4 | 5 |
|---|---|---|---|---|---|
| 1 | 1 | 4 | 2 | 5 | 3 |
| 2 | 4 | 2 | 5 | 3 | 1 |
| 3 | 2 | 5 | 3 | 1 | 4 |
| 4 | 5 | 3 | 1 | 4 | 2 |
| 5 | 3 | 1 | 4 | 2 | 5 |

**Example 8.3.** The following groupoid is 4-translatable, without idempotents.

|   | 1 | 2 | 3 | 4 | 5 |
|---|---|---|---|---|---|
| 1 | 2 | 1 | 3 | 4 | 5 |
| 2 | 1 | 3 | 4 | 5 | 2 |
| 3 | 3 | 4 | 5 | 2 | 1 |
| 4 | 4 | 5 | 2 | 1 | 3 |
| 5 | 5 | 2 | 1 | 3 | 4 |

**Example 8.4.** Re-examining Table two we see clearly that *a, ba, aba, ab, b* is a 3-translatable sequence for Q1= { *a, ab, ba, b, aba* }.

**Table 2.**

| Q1  | *a*  | *ab* | *ba* | *b*  | *aba* |        |        | (Q1)* | *a*    | *a*b*  | *b*a*  | *b*    | *aba*  |
|-----|------|------|------|------|-------|--------|--------|-------|--------|--------|--------|--------|--------|
| *a*   | *a*   | *ba*  | *aba* | *ab*  | *b*    | ........ | ........ | *a*     | *a*     | *b*     | *aba*   | *a*b*   | *b*a*   |
| *ab*  | *aba* | *ab*  | *b*   | *a*   | *ba*   | ........ | ........ | *a*b*   | *aba*   | *a*b*   | *b*     | *b*a*   | *a*     |
| *ba*  | *b*   | *a*   | *ba*  | *aba* | *ab*   | ........ | ........ | *b*a*   | *a*b*   | *a*     | *b*a*   | *aba*   | *b*     |
| *b*   | *ba*  | *aba* | *ab*  | *b*   | *a*    | ........ | ........ | *b*     | *b*a*   | *aba*   | *a*     | *b*     | *a*b*   |
| *aba* | *ab*  | *b*   | *a*   | *ba*  | *aba*  | ........ | ........ | *aba*   | *b*     | *b*a*   | *a*b*   | *a*     | *aba*   |

The sequence *a, a*b, aba, b*a, b* is a 2-translatable sequence for $(Q1)^* = \{\, a, b, b*a, aba, a*b\,\}$.

| (Q1)* | *a*  | *b*  | *b*a* | *aba* | *a*b* |
|-------|------|------|-------|-------|-------|
| *a*     | *a*   | *a*b* | *aba*  | *b*a*  | *b*    |
| *b*     | *b*a* | *b*   | *a*    | *a*b*  | *aba*  |
| *b*a*   | *a*b* | *aba* | *b*a*  | *b*    | *a*    |
| *aba*   | *b*   | *a*   | *a*b*  | *aba*  | *b*a*  |
| *a*b*   | *aba* | *b*a* | *b*    | *a*    | *a*b*  |

Note that from Table 2 we see that $(Q1)^*$ with ordering $(Q1)^* = \{\, a, a*b, b*a, b, aba\,\}$ is NOT translatable.

**Theorem 8.1.** *No idempotent, bookend, cancellable groupoid* $G = \{1, 2, ..., 9\}$ *of order 9 is translatable.*

Proof. If G has an *n*-translatable sequence $1, x_2, x_3, x_4, x_5, x_6, x_7, x_8, x_9$. If $n = 2$, then the first 4 Cayley rows are:

| G | 1 | 2 | 3 | 4 | 5 | 6 | 7 | 8 | 9 |
|---|---|---|---|---|---|---|---|---|---|
| 1 | 1 | $x_2$ | $x_3$ | $x_4$ | $x_5$ | $x_6$ | $x_7$ | $x_8$ | $x_9$ |
| 2 | $x_8$ | 2 | 1 | $x_2$ | $x_3$ | $x_4$ | $x_5$ | $x_6$ | $x_7$ |
| 3 | $x_6$ | $x_7$ | 3 | 2 | 1 | $x_2$ | $x_3$ | $x_4$ | $x_5$ |
| 4 | $x_4$ | $x_5$ | $x_6$ | 4 | 3 | 2 | 1 | $x_2$ | $x_3$ |

Then $1 \cdot 4 = x_4 = 4 \cdot 1 = x_4^2 = (1 \cdot 4)(4 \cdot 1) = 4 = (4 \cdot 1)(1 \cdot 4) = 1$, a contradiction.

If $n = 3$ then the first 4 rows of the Cayley table are:

|   | 1 | 2 | 3 | 4 | 5 | 6 | 7 | 8 | 9 |
|---|---|---|---|---|---|---|---|---|---|
| 1 | 1 | $x_2$ | $x_3$ | $x_4$ | $x_5$ | $x_6$ | $x_7$ | $x_8$ | $x_9$ |
| 2 | $x_7$ | 2 | $x_9$ | 1 | $x_2$ | $x_3$ | $x_4$ | $x_5$ | $x_6$ |
| 3 | $x_4$ | $x_5$ | 3 | $x_7$ | 2 | $x_9$ | 1 | $x_2$ | $x_3$ |
| 4 | 1 | $x_2$ | $x_3$ | . | . | . | . | . | . |

But then $1 = 1 \cdot 1 = 4 \cdot 1$ and so $1 = 4$, a contradiction.

If $n = 4$ then the first 4 rows of the Cayley table are:

|   | 1       | 2       | 3       | 4     | 5     | 6       | 7     | 8     | 9     |
|---|---------|---------|---------|-------|-------|---------|-------|-------|-------|
| 1 | 1       | $x_2$   | $x_3$   | $x_4$ | $x_5$ | $x_6$   | $x_7$ | $x_8$ | $x_9$ |
| 2 | $x_6$   | $x_7=2$ | $x_8$   | $x_9$ | 1     | $x_2$   | $x_3$ | $x_4$ | $x_5$ |
| 3 | $x_2$   | $x_3$   | $x_4=3$ | $x_5$ | $x_6$ | $x_7=2$ | $x_8$ | $x_9$ | 1     |
| 4 | $x_7=2$ | $x_8$   | $x_9$   | 1     | .     | .       | .     | .     | .     |

But then $4=4\cdot 4=1$, a contradiction.

If $n=5$ then the first 4 rows of the Cayley table are:

|   | 1     | 2     | 3     | 4     | 5     | 6     | 7     | 8     | 9     |
|---|-------|-------|-------|-------|-------|-------|-------|-------|-------|
| 1 | 1     | $x_2$ | $x_3$ | $x_4$ | $x_5$ | $x_6$ | $x_7$ | $x_8$ | $x_9$ |
| 2 | $x_5$ | $x_6$ | $x_7$ | $x_8$ | $x_9$ | 1     | $x_2$ | $x_3$ | $x_4$ |
| 3 | $x_9$ | 1     | $x_2$ | $x_3$ | $x_4$ | $x_5$ | $x_6$ | $x_7$ | $x_8$ |
| 4 | $x_4$ | $x_5$ | $x_6$ | $x_7$ | $x_8$ | .     | .     | .     | .     |

But then $4\cdot 1=1\cdot 4=x_4=x_4^2=(1\cdot 4)(4\cdot 1)=4=(4\cdot 1)(1\cdot 4)=1$, a contradiction.

If $n=6$ then the first 4 rows of the Cayley table are:

|   | 1     | 2     | 3     | 4     | 5     | 6     | 7     | 8     | 9     |
|---|-------|-------|-------|-------|-------|-------|-------|-------|-------|
| 1 | 1     | $x_2$ | $x_3$ | $x_4$ | $x_5$ | $x_6$ | $x_7$ | $x_8$ | $x_9$ |
| 2 | $x_4$ | $x_5$ | $x_6$ | $x_7$ | $x_8$ | $x_9$ | 1     | $x_2$ | $x_3$ |
| 3 | $x_7$ | $x_8$ | $x_9$ | 1     | $x_2$ | $x_3$ | $x_4$ | $x_5$ | $x_6$ |
| 4 | 1     | .     | .     | .     | .     | .     | .     | .     | .     |

But then $1=1\cdot 1=4\cdot 1$ and so $1=4$, a contradiction. If $n=7$ then the first 4 rows of the Cayley table are:

|   | 1     | 2     | 3     | 4     | 5     | 6     | 7     | 8     | 9     |
|---|-------|-------|-------|-------|-------|-------|-------|-------|-------|
| 1 | 1     | $x_2$ | $x_3$ | $x_4$ | $x_5$ | $x_6$ | $x_7$ | $x_8$ | $x_9$ |
| 2 | $x_3$ | $x_4$ | $x_5$ | $x_6$ | $x_7$ | $x_8$ | $x_9$ | 1     | $x_2$ |
| 3 | $x_5$ | $x_6$ | $x_7$ | $x_8$ | $x_9$ | 1     | $x_2$ | $x_3$ | $x_4$ |
| 4 | $x_7$ | $x_8$ | $x_9$ | 1     | .     | .     | .     | .     | .     |

But then $4=4\cdot 4=1$, a contradiction.

If $n=8$ then the first 2 rows of the Cayley table are:

|   | 1     | 2     | 3     | 4     | 5     | 6     | 7     | 8     | 9     |
|---|-------|-------|-------|-------|-------|-------|-------|-------|-------|
| 1 | 1     | $x_2$ | $x_3$ | $x_4$ | $x_5$ | $x_6$ | $x_7$ | $x_8$ | $x_9$ |
| 2 | $x_2$ | .     | .     | .     | .     | .     | .     | .     | .     |

But then $1\cdot 2=2\cdot 1=x_2=x_2^2=(1\cdot 2)(2\cdot 1)=2=(2\cdot 1)(1\cdot 2)=1$, a contradiction.

If $n=9$ then $2=2\cdot 2=x_2$ and $1\cdot 2=x_2=2$, which implies $1=2$, a contradiction.

By Proposition 8.1, G is not 1-translatable. Hence, G is not translatable. ∎

**Corollary 8.2.** *The quadratical quasigroups* $Q2$ *and* $(Q2)^*$ *are not translatable.*

**Theorem 8.3.** *The sequence* 11, 12, 33, 21, 31, 34, 24, 32, 13, 14, 23, *aba*, 22 *is a 5-translatable sequence for* $Q3 = \{11, 14, 34, 12, 23, 24, 33, aba, 32, 21, 22, 13, 31\}$. *The dual sequence*
$11^*, 12^*, 23^*, aba, 22^*, 13^*, 14^*, 34^*, 24^*, 32^*, 33^*, 21^*, 31^*$ *is an 8-translatable sequence for*
$(Q3)^* = \{11^*, 14^*, 31^*, 13^*, 21^*, 22^*, 33^*, aba, 32^*, 23^*, 24^*, 12^*, 34^*\}$.

Proof. Inspection of the Cayley table of Q3 (Table 3) yields the following table, which is clearly 5-translatable.

| Q3 | 11 | 14 | 34 | 12 | 23 | 24 | 33 | *aba* | 32 | 21 | 22 | 13 | 31 |
|---|---|---|---|---|---|---|---|---|---|---|---|---|---|
| 11 | **11** | **12** | **33** | **21** | **31** | **34** | **24** | **32** | **13** | **14** | **23** | *aba* | **22** |
| 14 | **13** | **14** | **23** | *aba* | **22** | **11** | **12** | **33** | **21** | **31** | **34** | **24** | **32** |
| 34 | **21** | **31** | **34** | **24** | **32** | **13** | **14** | **23** | *aba* | **22** | **11** | **12** | **33** |
| 12 | *aba* | **22** | **11** | **12** | **33** | **21** | **31** | **34** | **24** | **32** | **13** | **14** | **23** |
| 23 | **24** | **32** | **13** | **14** | **23** | *aba* | **22** | **11** | **12** | **33** | **21** | **31** | **34** |
| 24 | **12** | **33** | **21** | **31** | **34** | **24** | **32** | **13** | **14** | **23** | *aba* | **22** | **11** |
| 33 | **14** | **23** | *aba* | **22** | **11** | **12** | **33** | **21** | **31** | **34** | **24** | **32** | **13** |
| *aba* | **31** | **34** | **24** | **32** | **13** | **14** | **23** | *aba* | **22** | **11** | **12** | **33** | **21** |
| 32 | **22** | **11** | **12** | **33** | **21** | **31** | **34** | **24** | **32** | **13** | **14** | **23** | *aba* |
| 21 | **32** | **13** | **14** | **23** | *aba* | **22** | **11** | **12** | **33** | **21** | **31** | **34** | **24** |
| 22 | **33** | **21** | **31** | **34** | **24** | **32** | **13** | **14** | **23** | *aba* | **22** | **11** | **12** |
| 13 | **23** | *aba* | **22** | **11** | **12** | **33** | **21** | **31** | **34** | **24** | **32** | **13** | **14** |
| 31 | **34** | **24** | **32** | **13** | **14** | **23** | *aba* | **22** | **11** | **12** | **33** | **21** | **31** |

Inspection of the Cayley table of $(Q3)^*$ (Table 4) yields the required 8-translatable sequence for the groupoid $(Q3)^* = \{11^*, 14^*, 31^*, 13^*, 21^*, 22^*, 33^*, aba, 32^*, 23^*, 24^*, 12^*, 34^*\}$. ∎

**Theorem 8.4.** *The sequence* 11, 12, 42, 43, 13, 14, 33, 21, 31, 44, 23, *aba*, 22, 41, 34, 24, 32 *is a 13-translatable sequence for* $Q4 = \{11, 14, 23, 24, 43, 31, 41, 12, 33, aba, 32, 13, 44, 34, 42, 21, 22\}$ *and the dual sequence*
$11^*, 12^*, 34^*, 24^*, 32^*, 44^*, 23^*, aba, 22^*, 41^*, 33^*, 21^*, 31^*, 13^*, 14^*, 43^*, 42^*$ *is a 4-translatable sequence*
*for* $(Q4)^* = \{11^*, 14^*, 21^*, 22^*, 44^*, 34^*, 42^*, 13^*, 33^*, aba, 32^*, 12^*, 43^*, 31^*, 41^*, 23^*, 24^*\}$.

Proof. This follows from the inspection of the Cayley tables of $Q4$ and $(Q4)^*$, Tables 5 and 6 respectively. ∎

We now use the concept of a translatable groupoid to find a quadratical quasigroup of order 25, but not equal to $Q1 \times Q1$, $Q1 \times (Q1)^*$, $(Q1)^* \times Q1$ or $(Q1)^* \times (Q1)^*$. To avoid notational confusion, we will use italics to denote the positive integers *1, 2, ..., n*. So that *13* denotes the integer thirteen and NOT the element 13 = *ba* of a quadratical quasigroup.

**Theorem 8.5.** *No quadratical quasigroup of order 25 is 2, 3, 4, 5 or 6 translatable.*

Proof. Suppose that the sequence *1*, $x_2$, $x_3$, ..., $x_{25}$ is an *n*-translatable sequence for the quadratical quasigroup

$Q = \{1, 2, \ldots, 25\}$, with $n \in \{2, 3, 4, 5, 6\}$.

If $n = 2$ then $2 = (1 \cdot 2)(2 \cdot 1) = x_2 \cdot x_{24} = 25 \cdot 3 = x_5 = 22$, a contradiction.

If $n = 3$ then $2 = (1 \cdot 2)(2 \cdot 1) = x_2 \cdot x_{23} = 13 \cdot 15 = x_4 = 12$, a contradiction.

If $n = 4$ then $2 = (1 \cdot 2)(2 \cdot 1) = x_2 \cdot x_{22} = 9 \cdot 19 = x_{12} = 14$, a contradiction.

If $n = 5$ then $6 \cdot 1 = 1$, a contradiction.

If $n = 6$ then $6 \cdot 6 = 1$, a contradiction. ∎

**Theorem 8.6.** *Suppose that the Cayley Table of the idempotent groupoid* $K = \{1, 2, \ldots, 25\}$ *is determined by the 7-translatable sequence 1,5,9,13,17,21,25,4,8,12,16,20,24,3,7,11,15,19,23,2,6,10,14,18,22. Then K is a quadratical quasigroup with Cayley Table as follows, where $A = aba$:*

| K | 11 | 12 | 13 | 14 | A | 21 | 22 | 23 | 24 | 31 | 32 | 33 | 34 | 41 | 42 | 43 | 44 | 51 | 52 | 53 | 54 | x | xy | yx | y |
|---|---|---|---|---|---|---|---|---|---|---|---|---|---|---|---|---|---|---|---|---|---|---|---|---|---|
| 11 | 11 | 21 | A | 12 | 53 | 32 | yx | 54 | 51 | 44 | 34 | y | 52 | 23 | x | 14 | 31 | 24 | 41 | xy | 43 | 33 | 42 | 13 | 22 |
| 12 | A | 12 | 14 | 22 | 51 | 53 | 34 | 52 | x | yx | 43 | 54 | 33 | 13 | 21 | 32 | xy | y | 23 | 41 | 42 | 11 | 31 | 24 | 44 |
| 13 | 23 | 11 | 13 | A | 54 | y | 53 | 31 | 52 | 32 | 51 | 42 | xy | yx | 33 | 24 | 12 | 43 | 44 | 22 | x | 41 | 21 | 34 | 14 |
| 14 | 13 | A | 24 | 14 | 52 | 54 | 51 | xy | 33 | 53 | x | 31 | 41 | 34 | 11 | y | 22 | 42 | yx | 44 | 21 | 23 | 12 | 43 | 32 |
| A | 54 | 53 | 52 | 51 | A | 11 | 12 | 13 | 14 | 21 | 22 | 23 | 24 | 31 | 32 | 33 | 34 | 41 | 42 | 43 | 44 | y | yx | xy | x |
| 21 | 53 | yx | 51 | 34 | 12 | 21 | 31 | A | 22 | 42 | xy | 11 | 14 | 54 | 44 | x | 13 | 33 | y | 24 | 41 | 32 | 23 | 52 | 43 |
| 22 | 52 | 51 | 33 | x | 14 | A | 22 | 24 | 32 | 12 | 44 | 13 | y | xy | 53 | 11 | 43 | 23 | 31 | 42 | yx | 54 | 34 | 41 | 21 |
| 23 | y | 32 | 54 | 53 | 11 | 33 | 21 | 23 | A | x | 12 | 41 | 13 | 42 | 14 | 52 | yx | xy | 43 | 34 | 22 | 24 | 44 | 31 | 51 |
| 24 | 31 | 54 | xy | 52 | 13 | 23 | A | 34 | 24 | 11 | 14 | yx | 43 | 12 | y | 41 | 51 | 44 | 21 | x | 32 | 42 | 53 | 22 | 33 |
| 31 | 51 | 34 | x | 43 | 22 | 12 | xy | 14 | 44 | 31 | 41 | A | 32 | 52 | yx | 21 | 24 | 11 | 54 | y | 23 | 53 | 13 | 33 | 42 |
| 32 | xy | 52 | 41 | 33 | 24 | 13 | 14 | 43 | y | A | 32 | 34 | 42 | 22 | 54 | 23 | x | yx | 12 | 21 | 53 | 31 | 51 | 44 | 11 |
| 33 | 32 | 44 | 53 | yx | 21 | x | 42 | 11 | 12 | 43 | 31 | 33 | A | y | 22 | 51 | 23 | 52 | 24 | 13 | xy | 14 | 41 | 54 | 34 |
| 34 | 42 | y | 31 | 54 | 23 | 41 | 11 | yx | 13 | 33 | A | 44 | 34 | 21 | 24 | xy | 53 | 22 | x | 51 | 14 | 43 | 32 | 12 | 52 |
| 41 | 33 | x | 11 | 21 | 32 | 14 | 44 | y | 53 | 22 | yx | 24 | 54 | 41 | 51 | A | 42 | 13 | xy | 31 | 34 | 52 | 43 | 23 | 12 |
| 42 | 12 | 31 | 22 | xy | 34 | yx | 13 | 51 | 43 | 23 | 24 | 53 | x | A | 42 | 44 | 52 | 32 | 11 | 33 | y | 21 | 54 | 14 | 41 |
| 43 | yx | 23 | 34 | 13 | 31 | 42 | 4 | 12 | xy | y | 52 | 21 | 22 | 53 | 41 | 43 | A | x | 32 | 14 | 33 | 44 | 11 | 51 | 24 |
| 44 | 24 | 14 | y | 32 | 33 | 2 | x | 41 | 11 | 51 | 21 | xy | 23 | 43 | A | 54 | 44 | 31 | 34 | yx | 12 | 13 | 22 | 42 | 53 |
| 51 | 44 | 41 | yx | 23 | 42 | 43 | y | 21 | 31 | 24 | 54 | x | 12 | 32 | xy | 34 | 11 | 51 | 14 | A | 52 | 22 | 33 | 53 | 13 |
| 52 | x | 43 | 21 | 42 | 44 | 22 | 41 | 32 | yx | xy | 23 | 14 | 53 | 33 | 34 | 12 | y | A | 52 | 54 | 13 | 51 | 24 | 11 | 31 |
| 53 | 43 | 24 | 42 | y | 41 | xy | 33 | 44 | 23 | 52 | 11 | 22 | yx | x | 13 | 31 | 32 | 12 | 51 | 53 | A | 34 | 14 | 21 | 54 |
| 54 | 22 | xy | 44 | 41 | 43 | 34 | 24 | x | 42 | 13 | y | 51 | 21 | 14 | 31 | yx | 33 | 53 | A | 11 | 54 | 12 | 52 | 32 | 23 |
| x | 21 | 42 | 12 | 31 | yx | 44 | 23 | 53 | 34 | 41 | 13 | 32 | 51 | 11 | 43 | 22 | 54 | 14 | 33 | 52 | 24 | x | y | A | xy |
| xy | 14 | 22 | 32 | 44 | x | 51 | 43 | 33 | 21 | 34 | 42 | 52 | 11 | 24 | 12 | 53 | 41 | 54 | 13 | 23 | 31 | A | xy | y | yx |
| yx | 41 | 33 | 23 | 11 | y | 24 | 32 | 42 | 54 | 14 | 53 | 43 | 31 | 44 | 52 | 13 | 21 | 34 | 22 | 12 | 51 | xy | x | yx | A |
| y | 34 | 13 | 43 | 24 | xy | 31 | 52 | 22 | 41 | 54 | 33 | 12 | 44 | 51 | 23 | 42 | 14 | 21 | 53 | 32 | 11 | yx | A | x | y |

*Furthermore* $K \notin \{ Q1 \times Q1, Q1 \times (Q1)^*, (Q1)^* \times Q1, (Q1)^* \times (Q1)^* \}$

Proof. Consider the Cayley Table determined by the 7-translatable sequence *1, 5, 9, 13, 17, 21, 25, 4, 8, 12, 16, 20, 24, 3, 7, 11, 15, 19, 23, 2, 6, 10, 14, 18, 22* for the groupoid $K = \{1, 2, \ldots, 25\}$, where $a = 1, b = 2, x = 4$, $y = 24$, $xy = 9$ and $yx = 19$.

Note that the numbers in the table above are NOT written in italics and so, for example, $24 = b \cdot ba = 2 \cdot (2 \cdot 1) = 2 \cdot 23 = 11$. Similar calculations give the following values:

11 = *1*, 12 = *5*, 13 = *23*, 14 = *2*, *aba* = *14*, 21 = *17*, 22 = *18*, 23 = *10*, 24 = *11*, 31 = *21*, 32 = *15*, 33 = *13*, 34 = *7*, 41 = *22*, 42 = *8*, 43 = *20*, 44 = *6*, 51 = *16*, 52 = *25*, 53 = *3* and 54= *12*.

Table in Theorem 8.6 is then derived from these calculations. One still needs to check that the table is bookend and medial. These calculations are straightforward, but the number of them is great, so almost all are omitted. We only give the following few examples to demonstrate the methods of calculation.

For example, $(xy \cdot 23)(23 \cdot xy) = 33 \cdot 44 = 23$ and $(23 \cdot xy)(xy \cdot 23) = 44 \cdot 33 = xy$;
$(42 \cdot 54)(54 \cdot 42) = y \cdot 31 = 54$ and $(54 \cdot 42)(42 \cdot 54) = 31 \cdot y = 42$;
$(x \cdot 51)(13 \cdot yx) = 14 \cdot 34 = 41$ and $(x \cdot 13)(51 \cdot yx) = 12 \cdot 53 = 41$ and finally,
$(24 \cdot 43)(33 \cdot yx) = 41 \cdot 54 = 34$ and $(24 \cdot 33)(43 \cdot yx) = yx \cdot 51 = 34$.

We must now prove that $K \notin \{ Q1 \times Q1, Q1 \times (Q1)^*, (Q1)^* \times Q1, (Q1)^* \times (Q1)^* \}$

First we note that Q1 satisfies the identities $zw \cdot w = z$ and $z \cdot zw = wz$ and $(Q1)^*$ satisfies the dual identities $z \cdot zw = w$ and $zw \cdot w = wz$. Now, the subgroupoid of K generated by *x* and *y* is isomorphic to $(Q1)^*$.

If $K \cong Q1 \times Q1$, Q must satisfy the identity $z \cdot zw = wz$. Therefore, $x \cdot xy = yx$. But from the Cayley table of Theorem 8.6 we see that $x \cdot xy = y$. Therefore $y = yx$, a contradiction.

If $K \cong (Q1)^* \times (Q1)^*$ then K must satisfy the identity $z \cdot zw = w$. But from the Cayley table of Theorem 8.6 we see that $11 \cdot (11 \cdot 12) = 11 \cdot 21 = 32$, a contradiction.

Finally, $(Q1) \times (Q1)^*$ and $(Q1)^* \times (Q1)$ each contain five disjoint isomorphic copies of $(Q1)^*$. Now, $Z = \{x, xy, yx, y, A\}$, in K, is an isomorphic copy of $(Q1)^*$. However, from the Cayley table of Theorem 8.6 we see that no other two elements of $K - Z$ generate a subgroupoid of K. Hence, K is NOT isomorphic to $(Q1)^* \times (Q1)$ or $(Q1) \times (Q1)^*$. ∎

**Corollary 8.7.** *The dual groupoid* $K^* = \{1,2, \ldots ,25\}$ *of the quadraical quasigroup* K *of Theorem 8.6 has an 18-translatable sequence 1,23,20,17,14,11,8,5,2,24,21,18,15,12,9,6,3,25,22,19,16,13,10,7,4.*

Proof. This follows from the fact that the entries of the *k*-th row of the dual groupoid, from left to right, are exactly the entries of the *k*-th column of K, from top to bottom, where $k \in \{1, 2, ..., 25\}$. ∎

The proof of the following Theorem is omitted, for reasons of length.

**Theorem 8.8.** *The sequence*
11,12,54,74,43,62,53,44,23,*aba*,22,41,52,63,42,71,51,13,14,61,33,21,31,73,72,34,24,32,64 *is a 12-translatable sequence for the quadratical quasigroup*
$G_{29}$ = {11,14,44,34,42,74,21,22,64,13,53,72,63,54,33,*aba*,32,51,62,73,52,12,61,23,24,71,43,31,41} *of form Q7.*
*The groupoid* $(G_{29})^*$ *is 17-translatable.*

## 9. Open questions and directions for future research.

The following table summarises some of our results.

| Groupoid | Order | Of form $Qn$? | $k$-translatable? | Generated by any two distinct elements? | 2-generated? | Self-dual? |
|---|---|---|---|---|---|---|
| Q1 | 4 | $n=1$ | $k=3$ | Yes | Yes | No |
| $(Q1)^*$ | 4 | $n=1$ | $k=2$ | Yes | Yes | No |
| $Q2 \cong (Q2)^*$ | 9 | $n=2$ | No | Yes | Yes | Yes |
| Q3 | 13 | $n=3$ | $k=5$ | Yes | Yes | No |
| $(Q3)^*$ | 13 | $n=3$ | $k=8$ | Yes | Yes | No |
| Q4 | 17 | $n=4$ | $k=13$ | Yes | Yes | No |
| $(Q4)^*$ | 17 | $n=4$ | $k=4$ | Yes | Yes | No |
| $Q1 \times Q1$ | 25 | No | No | No | No | No |
| $(Q1)^* \times (Q1)^*$ | 25 | No | No | No | No | No |
| $(Q1) \times (Q1)^*$ | 25 | No | No | No | Yes | No |
| $(Q1)^* \times (Q1)$ | 25 | No | No | No | Yes | No |
| K | 25 | No | $k=7$ | No | Yes | No |
| $K^*$ | 25 | No | $k=18$ | No | Yes | No |
| $G_{29}$ | 29 | $n=7$ | $k=12$ | Yes | Yes | No |
| $(G_{29})^*$ | 29 | $n=7$ | $k=17$ | Yes | Yes | No |

**Open questions:**

1. For which values of $n$ is there a quadratical groupoid of form $Qn$ ? (Note that $n \notin \{5,6\}$ ).
2. Is every quadratical groupoid of form $Qn$ $(n \geq 7)$ translatable?
3. Are there self-dual, quadratical groupoids of order greater than 9?
4. Is every quadratical groupoid of order greater than 9 and of form $Qn$ $(n \geq 3)$ generated by any two of its distinct elements?
5. If a quadratical quasigroup Q of order $n$ is $k$-translatable then is $Q^*$ $(n-k)$-translatable?
6. For which values of $n$ and $k$ is there a $k$-translatable, quadratical groupoid of order $n$?
7. Is there a self-dual, quadratical and translatable quasigroup?
8. Conjecture: Does $Sp\ Q = \{\ 1+4n\ :\ n \in \mathbb{N}_0 - \{5\}\ \}$ ?
9. If an idempotent groupoid $G = \{1, 2, ..., n\}$ has a $k$-translatable sequence $1, x_2, x_3, ..., x_{n-1}, x_n$ and $2 = (1 \cdot 2)(2 \cdot 1)$ then is G bookend? Is it a quadratical quasigroup?
10. Do translatable, idempotent semigroups exist?
11. Determine **all** the remaining inter-relationships amongst the properties of idempotency, elasticity, strong elasticity, bookend, left-distributivity, right-distributivity, left-cancellativity, right-cancellativity, mediality, alterability, right-solvability, left-solvability and Property A.

[1]Flat 10, Albert Mansions, Crouch Hill, London N8 9RE, United Kingdom
bobmonzo@talktalk.net